# CONCENTRATION OF MEASURE AND SPECTRA OF RANDOM MATRICES: APPLICATIONS TO CORRELATION MATRICES, ELLIPTICAL DISTRIBUTIONS AND BEYOND[1]

By Noureddine El Karoui

*University of California, Berkeley*

We place ourselves in the setting of high-dimensional statistical inference, where the number of variables $p$ in a data set of interest is of the same order of magnitude as the number of observations $n$. More formally, we study the asymptotic properties of correlation and covariance matrices, in the setting where $p/n \to \rho \in (0, \infty)$, for general population covariance.

We show that, for a large class of models studied in random matrix theory, spectral properties of large-dimensional correlation matrices are similar to those of large-dimensional covariance matrices.

We also derive a Marčenko–Pastur-type system of equations for the limiting spectral distribution of covariance matrices computed from data with elliptical distributions and generalizations of this family. The motivation for this study comes partly from the possible relevance of such distributional assumptions to problems in econometrics and portfolio optimization, as well as robustness questions for certain classical random matrix results.

A mathematical theme of the paper is the important use we make of concentration inequalities.

**1. Introduction.** It is increasingly common in multivariate statistics and various areas of applied mathematics and computer science to have to work with data sets where the number of variables, $p$, is of the same order of magnitude as the number of observations, $n$. When studying asymptotic properties of estimators in this setting, usually under the assumption that $p/n$ has a finite nonzero limit, we often obtain convergence results that differ

Received October 2007; revised May 2008.

[1]Supported by NSF Grant DMS-06-05169 and SAMSI.

from those obtained under the "classical" assumptions that $p$ is fixed and $n$ goes to infinity.

A good example is provided by problems of portfolio optimization in quantitative finance. Typically, if one is working with, say, the stocks making up the S&P 500 index over a period of one year and recording quantities daily, one will have to work with a data matrix of size roughly $250 \times 500$. For many "large" portfolio optimization problems, the data matrices will have the characteristic that $p/n$ is not very small. Since covariance matrices (and their inverses) play a key role in solving a number of portfolio optimization problems and, in particular, Markowitz's formulation (see, e.g., [13]), it is important to understand how well our estimators perform in this "new" type of asymptotics. This should help us to assess the quality of our empirical choices of portfolio and their proximity to the theoretically optimal ones. We refer the reader who is particularly interested in these questions to [31] for an early and interesting application of random matrix ideas to portfolio optimization problems.

The realization of the fact that there might be problems in estimating the spectral properties of large-dimensional covariance matrices when $p/n$ is not small is not recent: the first paper in the area is probably [35], where the authors studied the behavior of the eigenvalues of large-dimensional sample covariance matrices for diagonal population covariance matrices and with some assumptions on the structure of the data. The surprising result they found was that, in the case of i.i.d. data with variance 1, the eigenvalues of the sample covariance matrix $X^*X/n$ do not concentrate around 1 (the value of all population eigenvalues), but rather were spread out on the interval $[(1 - \sqrt{p/n})^2, (1 + \sqrt{p/n})^2]$ when $p \leq n$. Moreover, their empirical distribution is asymptotically nonrandom. We note that this seminal paper is much richer than just described and refer the reader to it for more details. A simple lesson to be taken from this is that when $p/n$ is not small, the sample covariance matrix is not a good estimator of the population covariance.

Since this result, there has been a flurry of activity, especially in recent years, concerning the behavior of the largest eigenvalue of sample covariance matrices [21, 50], their fluctuation behavior in the null case [14, 19, 28, 29] and under alternatives [6, 16, 40], as well as fluctuation results for linear spectral statistics of those matrices [1, 5, 30]. Even more recently, some of these results have started to be used to develop better estimators of these large-dimensional covariance matrices ([11, 17] and [42]). We also note that from a statistical point of view, other approaches to estimation using regularization have been taken, with sometime striking results [8, 32].

As noted above, the random matrix results in question concern, somewhat exclusively, sample covariance matrices. However, in practice, sample correlation matrices are often used, for instance for principal component analysis (PCA). A question we were asked several times by practitioners is the extent to which the random matrix results would hold if one were concerned with correlation matrices and not covariance matrices. Part of the answer is already known from a paper due to [27], where he considered the case of i.i.d. data. The answer was that spectral distribution results, as well as a.s. convergence of extreme eigenvalue results, held in this situation. However, in practice, the assumption of i.i.d. data is not very reasonable and, in most cases, practitioners would actually hope to be in the presence of an interesting covariance structure, away from the no-information case represented by the identity covariance matrix. In this paper, we tackle the case where the population covariance is not $\mathrm{Id}_p$ and show that classic random matrix results hold then too, with the population covariance matrix replaced by the population correlation matrix. This means that recently developed methods that make use of random matrix theory to better estimate the eigenvalues of population covariance matrices can also be used to estimate the spectrum of population correlation matrices.

As explained below, such results can be shown for Gaussian and some non-Gaussian data. Therefore, a natural question is to wonder how robust to these distributional assumptions the results are. In particular, a recent paper [20] and a recent monograph [37] make an interesting case for modeling financial data through elliptical distributions. As explained in [20] and [37], this has to do with certain tail-dependence properties that are absent from Gaussian data and present in a certain class of elliptically distributed data. As mentioned above, understanding the spectral properties of sample covariance matrices with these distributions should help us better understand the properties of empirical solutions to the classical Markowitz portfolio optimization problem. This is one of the many applications these results could have. Though this paper does not deal with this specific problem, in the second part of the paper, we show that for elliptically distributed data (and generalizations), the spectrum of the sample covariance matrix is asymptotically nonrandom and we characterize the limit through the use of Stieltjes transforms. In particular, the result shows that the Marčenko–Pastur equation is not robust to deviation from the "Gaussian+" model usually considered in random matrix theory (see [44] and Theorem 1 below for an example of those assumptions). The result also explains some of the numerical results obtained by [20]. From a more theoretical standpoint, our approach allows us to break away from models for which the data vectors are linear transformations of random vectors with independent entries. Rather, what we need are concentration properties for 1-Lipschitz (with respect to the Euclidean norm) functionals of these data vectors. We note that some of our results can be obtained when the concentration properties are limited to convex 1-Lipschitz functions. Hence, our approach will show that some classical results in random matrix theory hold in wider generality than was previously known. For instance, it shows that classical random matrix results apply to data drawn, for instance, from a Gaussian copula (see [38]), under some restrictions on the operator norm of the corresponding correlation matrix. We note that the Gaussian copula problem appears to be quite far from what can be obtained using currently available results.

As it turns out, central to the proofs to be presented are the concentration properties of certain quadratic forms. Below, we make use of a number of concentration inequalities, recent and less recent. The usefulness of these inequalities in random matrix theory has already been illustrated in [26], in a different context from that which we develop below. A very good reference on the topic of concentration is [33].

The fact that we rely on concentration of quadratic forms for many of these results also yields some practical insights about possible limitations of the models considered in the random matrix literature. In particular, applying the concentration results to the standard random matrix models considered in, for example, Theorem 1 (or their generalizations in Theorem 2, with $\lambda_i = 1$ in the notation of this latter theorem) shows that the corresponding data vectors have norms (when divided by the square root of the dimension) that are all almost equal. Similarly, one can show (see, e.g., [15] for more details) that the concentration results we need also imply that for standard models, the data vectors are almost orthogonal to one another. More precisely, one can show that the maximum angle between data vectors goes to zero almost surely. Before applying or using these random matrix results and, in particular, the Marčenko–Pastur equation, it therefore appears that practitioners should pay attention to these features of the data by, for instance, drawing histograms of the angles between data vectors and norms of the data vector divided by the square root of the dimension. If those are not "concentrated," this might call into question the quality of the fit of the standard models and the relevance of the insights drawn from random matrix results.

The paper is organized as follows. In Section 2, we study the problem of spectral characteristics of correlation matrices with data drawn from standard random matrix models. In Section 3, we characterize the limiting spectrum of covariance matrices computed from data drawn from a generalization of elliptical distributions. In terms of concentration results, Section 2 can be viewed as using concentration results for the norm of the columns of the data matrices of interest. On the other hand, Section 3 relies on concentration properties for the rows of the data matrices of interest.

**2. On large dimensional correlation matrices.** We recall that a correlation matrix is a matrix that contains the correlations between the entries of a vector. So, if $R$ is the correlation matrix of the vector $v$,

$$R_{i,j} = \frac{\operatorname{cov}(v_i, v_j)}{\sqrt{\operatorname{var}(v_i)}\sqrt{\operatorname{var}(v_j)}}.$$

We will naturally focus on empirical correlation matrices. We assume that we are given a sample of data $\{y_k\}_{k=1}^n$, where $y_k \in \mathbb{R}^p$. Let us call $\bar{r}_{\cdot,i} = 1/n \sum_{k=1}^n y_{k,i}$ the mean of the $i$th component of our data vectors. We call $s_i$ the standard estimate of standard deviation of $\{y_{k,i}\}_{k=1}^n$, that is,

$$s_i^2 = \frac{1}{n-1} \sum_{i=1}^n (y_{k,i} - \bar{y}_{\cdot,i})^2.$$

By definition, $\widehat{R}$, the empirical correlation matrix of the data, is

$$\widehat{R}_{i,j} = \frac{1/(n-1) \sum_{k=1}^n (y_{k,i} - \bar{y}_{\cdot,i})(y_{k,j} - \bar{y}_{\cdot,j})}{s_i s_j}.$$

These matrices play an important role in many multivariate statistical methods. In particular, in techniques like principal component analysis, there are sometimes debates as to whether one should use the correlation matrix of the data or their covariance matrix. It is therefore important for practitioners to have information about the behavior of correlation matrices in high dimensions.

We now turn to our study of sample correlation matrices. The main result is Theorem 1, which states that under the model considered there (related to the classical one in random matrix theory), results concerning the spectral distribution and the largest eigenvalue carry over, without much modification, from sample covariance matrices to sample correlation matrices.

Before we proceed, we need to establish some notation. In the remainder of the paper, $\mathbb{C}^+ = \{z \in \mathbb{C} : \operatorname{Im}[z] > 0\}$. We call $v'$ the *transpose* of the vector $v$ and use the same notation for matrices. We use $\|M\|_2$ to denote the operator norm of a matrix $M$, that is, its largest singular value. For a positive semidefinite matrices, it is obviously its largest eigenvalue. If $Y$ is an $n \times p$ matrix, we naturally denote by $Y_{i,j}$ its $(i, j)$ entry and call $\bar{Y}$ the matrix whose $j$th column is constant and equal to $\bar{Y}_{\cdot,j}$. Finally, the sample covariance matrix of the data stored in matrix $Y$ is

$$S_p = \frac{1}{n-1}(Y - \bar{Y})'(Y - \bar{Y}).$$

2.1. *A simple lemma.* The crux of our argument is going to be that correlation matrices can be represented as the products of certain matrices, one of them being a covariance matrix. Hence, we will have solved our problem if we can show that these other matrices are "not too far" from matrices we understand well and if we can show that "nothing" (as far as spectral properties are concerned) is lost when replacing them by these better understood matrices. Before we state our main theorem and prove it, we state two results of independent interest, on which we will rely in the proof.

LEMMA 1. *Suppose that $M_p$ is a $p \times p$ Hermitian random matrix whose spectral characteristics [spectral distribution $F_p$ or largest eigenvalue $\lambda_1(M_p)$] converge a.s. to a limit and whose spectral norm is (a.s.) bounded as $p \to \infty$. Suppose that $D_p$ is a $p \times p$ diagonal matrix and that $\|\|D_p - \mathrm{Id}_p\|\|_2 \to 0$ a.s. Then the spectral characteristics of $D_p M_p D_p$ and $D_p^{-1} M_p D_p^{-1}$ have the same limits as those of $M_p$.*

PROOF. The assumption $\|\|D_p - \mathrm{Id}_p\|\|_2 \to 0$ implies that for $p$ large enough, $D_p$ is invertible. Now,

$$\|\|M_p - D_p M_p D_p\|\|_2 = \|\|M_p - M_p D_p + M_p D_p - D_p M_p D_p\|\|_2$$
$$\leq \|\|M_p\|\|_2 \|\|D_p - \mathrm{Id}_p\|\|_2 + \|\|D_p - \mathrm{Id}_p\|\|_2 \|\|D_p\|\|_2 \|\|M_p\|\|_2$$
$$\to 0 \qquad \text{a.s.}$$

Using Weyl's inequality (see [7], Corollary III.2.6), that is, the fact that for Hermitian matrices $A$ and $B$, and any $i$, if $\lambda_i(A)$ denotes the $i$th eigenvalue of $A$, ordered in decreasing order, $|\lambda_i(A) - \lambda_i(B)| \leq \|\|A - B\|\|_2$, we conclude that

$$\max_{k=1,\dots,p} |\lambda_k(M_p) - \lambda_k(D_p M_p D_p)| \to 0 \qquad \text{a.s.}$$

Because $\|\|M_p\|\|_2$ is bounded a.s., the two sequences are a.s. asymptotically distributed (see [25], page 62, or [24]). Therefore, if $F_p(M_p)$ converges weakly to $F$, then $F_p(D_p M_p D_p)$ also converges to $F$.

If $\|\|D_p - \mathrm{Id}_p\|\|_2 \to 0$, then $\|\|D_p^{-1} - \mathrm{Id}_p\|\|_2 \to 0$ too. So, the same results hold when we replace $D_p$ by $D_p^{-1}$. □

The previous lemma is helpful in our context thanks to the following elementary fact, which is standard in multivariate statistics.

FACT 1 (Correlation matrix as function of covariance matrix). *Let $C_p$ denote the correlation matrix of our data and $S_p$ the covariance matrix of the data. Let $D_p(S_p)$ denote the diagonal matrix consisting of the diagonal of $S_p$. We then have*

$$C_p = [D_p(S_p)]^{-1/2} S_p [D_p(S_p)]^{-1/2}.$$

PROOF. This is just a simple consequence of the fact that if $D$ is a diagonal matrix, then

$$(DHD)_{i,j} = d_{i,i} H_{i,j} d_{j,j}.$$

Note that $C_p(i,j) = S_p(i,j)/\sqrt{S_p(i,i)S_p(j,j)}$ and the assertion follows. □

As a consequence of the previous lemma and fact, we will deduce the asymptotic spectral properties of correlation matrices from those of covariance matrices by simply showing convergence of the diagonal of $S_p$ (or a scaled version of it) to $\mathrm{Id}_p$ in operator norm.

2.2. *Spectra of large-dimensional correlation matrices.* We are now ready to state the main theorem of this section.

THEOREM 1. *Suppose that $X$ is an $n \times p$ matrix of i.i.d. random variables with variance 1. Assume, without loss of generality, that their common mean is 0. Denote by $X_{i,j}$ the $(i,j)$th entry of $X$. Assume, further, that $\mathbf{E}(|X_{i,j}|^4(\log(|X_{i,j}|))^{2+2\varepsilon}) < \infty$. Suppose that $\Sigma_p$ is a $p \times p$ covariance matrix and let $\Gamma_p$ denote the corresponding correlation matrix. Assume that $\|\Gamma_p\|_2 < K$ for all $p$.*

*Let*

- *$Y = X\Sigma_p^{1/2}$ ($Y$ is the observed data matrix—the $n$ observed data vectors are stored in the rows of $Y$);*
- *$Y_1 = X\Gamma_p^{1/2}$.*

*The spectral properties of $\mathrm{corr}(Y)$, the sample correlation matrix of the data, are then the same as the spectral properties of $\Gamma_p^{1/2}(X - \bar{X})'(X - \bar{X})\Gamma_p^{1/2}/(n-1) = (Y_1 - \bar{Y}_1)'(Y_1 - \bar{Y}_1)/(n-1)$.*

*In particular, the Stieltjes transform of the limiting spectral distribution of $\mathrm{corr}(Y)$ satisfies the Marčenko–Pastur equation, with parameter the spectral distribution of $\Gamma_p$. Namely, if $H_p$, the spectral distribution of $\Gamma_p$, has a.s. a limit $H$, if $p/n$ has a finite limit $\rho$ and if $m_n$ is the Stieltjes transform of $\mathrm{corr}(Y)$, we have, letting $w_n = -(1-p/n)/z + (p/n)m_n(z)$,*

$$w_n(z) \to w(z) \qquad a.s., \text{ which satisfies } -\frac{1}{w(z)} = z - \rho \int \frac{\lambda \, dH(\lambda)}{1 + \lambda w(z)},$$

*and $w$ is the unique function mapping $\mathbb{C}^+$ into $\mathbb{C}^+$ to satisfy this equation.*

*Also, if the norm of $\Gamma_p^{1/2}(X - \bar{X})'(X - \bar{X})\Gamma_p^{1/2}/(n-1)$ has a limit in which $\Gamma_p$ intervenes only through its eigenvalues, the norm of $\mathrm{corr}(Y)$ has the same limit.*

This theorem is related to that of [27], which was concerned with $\Gamma_p = \mathrm{Id}_p$, which would amount to doing multivariate analysis with i.i.d. variables, an assumption that, for obvious statistical reasons, practitioners are not willing to make. Here, by contrast, we are able to handle general covariance structures, assuming that the spectral norm of $\Gamma_p$ is bounded. However, [27] required only four moments and we require a little more. We explain in Section 2.3.2 why this is the case.

We note that the proof can actually handle cases where $\|\Gamma_p\|_2$ grows slowly with $p$. We refer the reader to [44] for more information on the Marčenko–Pastur equation. We note that the paper [44] is an important strengthening of the result of [35], dealing, in particular, with nondiagonal covariance matrices.

Recent progress has led to fairly explicit characterization of the norm of large-dimensional sample covariance matrices, a fact that makes these results potentially useful in, among other fields, statistics. In particular, results concerning limiting spectral distributions do not, in general, provide any information about the localization of the largest eigenvalue of the corresponding matrices. For many practitioners and, in particular, those dealing with principal component analysis (see, e.g., [36]), it is important to have this localization information. Our analysis, combined with recent results, allows us to characterize the limit of the largest eigenvalue of sample correlation matrices in certain cases.

In particular, the following consequence for the norm of the correlation matrix can be drawn from the recent article [16], specifically Fact 2 there (which is partly a consequence of a deep result in [4]).

COROLLARY 1. *Under the assumptions of Theorem 1, if $\lambda_1(\Gamma_p)$ tends to the endpoint of the support of $H$ and the model $\{\Gamma_p, n, p\}$ is in the class $\mathcal{G}$ defined in [16], then*

$$\|\operatorname{corr}(Y)\|_2 - \mu_{n,p} \to 0 \qquad a.s.,$$

*where*

$$\mu_{n,p} = \frac{1}{c_0}\left(1 + \frac{p}{n}\int \frac{\lambda c_0}{1 - \lambda c_0}\, dH(\lambda)\right),$$

$$\frac{n}{p} = \int \left(\frac{\lambda c_0}{1 - \lambda c_0}\right)^2 dH(\lambda), \qquad c_0 \in [0, 1/\lambda_1(\Gamma_p)).$$

Although the result might seem somewhat cryptic (the conditions set forth in [16] are rather complicated to describe), it basically says that if the largest eigenvalues of $\Gamma_p$ are sufficiently close to one another (the precise mathematical meaning of this statement is contained in the assumptions made in [16]), then the largest eigenvalue of $C_p$ will converge to the endpoint of the limiting spectral distribution of $\text{corr}(Y)$. We insist on the fact that the quantities above ($c_0$ and $\mu_{n,p}$) are fairly easy to compute explicitly with a computer, making them relevant in practice. We refer the reader to [16] for more information about these problems, as well as examples of matrices $\Gamma_p$ for which the results hold.

2.3. *Proof of Theorem 1.* The proof is in three steps. The first involves showing that to understand the spectral properties of $\text{corr}(Y)$, we need to focus only on the matrix $Y_1'Y_1$ [or $(Y_1 - \bar{Y}_1)'(Y_1 - \bar{Y}_1)$]. We then need a truncation and centralization step for the entries of $X$. Finally, we use a concentration-of-measure results to show that the diagonal of the corresponding covariance matrix indeed converges in operator norm to the identity. We postpone a formal proof of Theorem 1 to Section 2.3.5, where we put all of the elements together.

2.3.1. *Replacing $\Sigma_p$ by $\Gamma_p$.* Since the correlation coefficient is invariant under shifting and (positive) scaling of random variables, we see that for any diagonal matrix $D$ with positive entries,

$$\text{corr}(Y) = \text{corr}(YD)$$

since $(YD)_{i,j} = Y_{i,j}d_{jj}$. In particular, for $D$, we can use $(\text{diag}(\Sigma_p))^{-1/2}$, which clearly has positive entries. After this adjustment, the data matrix we will focus on, $Y_2$, takes the form

$$Y_2 = XG, \qquad \text{where } G = \Sigma_p^{1/2}(\text{diag}(\Sigma_p))^{-1/2} = \Sigma_p^{1/2}D,$$

and $G'G = \Gamma_p$. Note, in particular, that since $\Gamma_p$ is a correlation matrix, its diagonal consists of 1's. Because $G$ is not symmetric, it is not, in general, equal to $\sqrt{\Gamma_p}$. We thus need to explain why we will be able to rely on existing random matrix results since, for example, [44] requires the data to have the form $X\Sigma_p^{1/2}$.

Since $G$ is similar to $D^{1/2}\Sigma_p^{1/2}D^{1/2}$, all of its eigenvalues are real and nonnegative. Further, because $G'G = \Gamma_p$, the eigenvalues of $G$ are equal to the square root of the eigenvalues of $\Gamma_p$. Because $\Sigma_p^{1/2}$ and $D$ are invertible, so is $\Sigma_p^{1/2}D$. Therefore, the spectrum of the matrix of interest, $Y_2'Y_2/n$, is the same as the spectrum of $X'X\Sigma_p^{1/2}D^2\Sigma_p^{1/2}/n$. Even though, in general, $\Sigma_p^{1/2}D^2\Sigma_p^{1/2} \neq \Gamma_p$, these matrices have the same eigenvalues. Because the Marčenko–Pastur equation involves only the eigenvalues of the deterministic matrix in question, the limiting spectral distribution of $Y_2'Y_2/n$ is the same as the limiting spectral distribution of $\Gamma_p^{1/2}X'X\Gamma_p^{1/2}/n = Y_1'Y_1/n$. A similar conclusion applies to the largest eigenvalue if it depends only of $\Gamma_p$ through $\Gamma_p$'s spectrum.

So, in what follows, we only need to investigate $\text{corr}(Y_2)$ or $\text{corr}(Y_1)$ to understand $\text{corr}(Y)$.

2.3.2. *Truncation and centralization step.* In this subsubsection, we show that we can truncate the entries of $X$, the $n \times p$ matrix full of i.i.d. random variables, at level $\sqrt{n}/(\log n)^{(1+\varepsilon)/2} = \sqrt{n}/\delta_n$ and almost surely not change the value of corr($Y$), at least for $p$ large enough. The same holds when the truncated values are then re-centered. The conclusion of this subsubsection is that it is enough to study matrices $X$ whose entries are i.i.d. with mean 0 and bounded in absolute value by $C\sqrt{n}/(\log n)^{(1+\varepsilon)/2}$.

The proof is similar to the argument given for the proof of Lemma 2.2 in [50]. However, because the term $1/(\log n)^{(1+\varepsilon)/2}$ is crucial in our later arguments and the authors of [50] gloss over the details of their choice of $\delta_n$, we feel a full argument is needed to give a convincing proof, although we do not claim that the arguments are new. This is where we need a slightly stronger assumption that just the finite fourth moment assumption made in [50]. (Our problem is with Remark 1 in [50], which is not clearly justified. There also appears to be counterexamples to this claim. However, it does not seem that (the full strength of) this remark is ever really used in that paper and the rest of the arguments are clear.) We have the following lemma, which closely follows Lemma 2.2 in [50].

LEMMA 2 (Truncation). *Let $X$ be an infinite double array of identically distributed (i.d.) random variables. Suppose that $X_n$ is an $n \times p$ matrix of identically distributed random variables, with mean 0, variance 1 and whose entries, $X_{i,j}$, satisfy $\mathbf{E}(|X_{i,j}|^4 (\log(|X_{i,j}|))^{2+2\varepsilon}) < \infty$. $X_n$ corresponds to the upper-left corner of $X$. Suppose that $p/n$ has a finite limit $\rho$. Let $T_n$ denote the matrix with $(i, j)$th entry $X_{i,j} 1_{|X_{i,j}| < \sqrt{n}/(\log n)^{(1+\varepsilon)/2}}$. Then,*

$$P(X_n \neq T_n \ i.o.) = 0.$$

PROOF. Because of the moment assumption made on $X_{i,j}$, we have, if we let $f_\varepsilon(x) = x^4 (\log x)^{2(1+\varepsilon)}$,

$$\int_0^\infty f'_\varepsilon(y) P(|X_{i,j}| > y) \, dy = \sum_{m=0}^\infty \int_{u_m}^{u_{m+1}} f'_\varepsilon(y) P(|X_{i,j}| > y) \, dy < \infty$$

for any increasing sequence $\{u_m\}_{m=0}^\infty$, with $u_0 = 0$ and $u_m \to \infty$ as $m \to \infty$. Now, when $y$ is large enough, $f'_\varepsilon(y) \geq 0$, so

$$\int_{u_m}^{u_{m+1}} f'_\varepsilon(y) P(|X_{i,j}| > y) \, dy \geq P(|X_{i,j}| > u_{m+1})(f_\varepsilon(u_{m+1}) - f_\varepsilon(u_m)).$$

Let $\gamma_m = 2^m$ and $u_m = \sqrt{\gamma_m/(\log \gamma_m)^{1+\varepsilon}}$. Note that $u_m$ is increasing for $m$ sufficiently large. Elementary computations show that as $m$ tends to $\infty$, $u_m^4 (\log u_m)^{2+2\varepsilon} \sim 2^{2m-(2+2\varepsilon)}$. Consequently, $f_\varepsilon(u_{m+1}) - f_\varepsilon(u_m) \sim 3 \times 2^{2(m-1)}$. Note that our moment requirements therefore imply that

$$\sum_{m=1}^\infty 2^{2m} P(|X_{i,j}| > u_m) < \infty.$$

Now, for $n$ satisfying $\gamma_{m-1} \le n < \gamma_m$, we threshold $X_n(i,j)$ at level $u_{m-1}$. (In what follows, $2\rho\gamma_m$ should be replaced by the smallest integer greater than this number, but to avoid cumbersome notation, we do not stress this particular fact.) We have

$$P(X_n \neq T_n \text{ i.o.}) \le \sum_{m=k}^{\infty} P\left( \bigcup_{\gamma_{m-1} \le n < \gamma_m} \bigcup_{i=1}^{n} \bigcup_{j=1}^{p} (|X_n(i,j)| > u_{m-1}) \right)$$

$$\le \sum_{m=k}^{\infty} P\left( \bigcup_{\gamma_{m-1} \le n < \gamma_m} \bigcup_{i=1}^{\gamma_m} \bigcup_{j=1}^{2\rho\gamma_m} (|X_n(i,j)| > u_{m-1}) \right)$$

$$= \sum_{m=k}^{\infty} P\left( \bigcup_{i=1}^{\gamma_m} \bigcup_{j=1}^{2\rho\gamma_m} (|X_n(i,j)| > u_{m-1}) \right)$$

$$\le 2\rho \sum_{m=k}^{\infty} \gamma_m^2 P(|X_{i,j}| > u_{m-1})$$

$$= 8\rho \sum_{m=k}^{\infty} 2^{2(m-1)} P(|X_{i,j}| > u_{m-1}).$$

The right-hand side tends to 0 when $k$ tends to infinity and the left-hand side is independent of $k$. We conclude that

$$P(X_n \neq T_n \text{ i.o.}) = 0. \qquad \square$$

LEMMA 3 (Centralization). *Let $TC_n$ denote the matrix with entries $TC_n(i, j) = T_n(i,j) - ET_n(i,j)$. Then,*

$$\frac{1}{n} \|T_n' T_n - TC_n' TC_n\|_2 \to 0 \qquad a.s.$$

PROOF. The proof would be a simple repetition of the arguments in the proof of Lemma 2.3 in [50], with $r = 1/2$ and $\delta = (\log n)^{-(1+\varepsilon)/2}$ in the notation of their papers, so we omit it. Note that that proof finds a bound on the spectral norm of $T_n' T_n - TC_n' TC_n$. $\square$

The centralization lemma (Lemma 3) guarantees that the spectral characteristics of $G' TC_n' TC_n G/n$ are asymptotically the same as those of $G' T_n' T_n G/n$: this is a consequence of the fact that $\|G\|_2^2 = \lambda_1(\Gamma_p)$ is uniformly bounded, as well as of Weyl's inequality,

$$\max_i |\lambda_i(G' TC_n' TC_n G/n) - \lambda_i(G' T_n' T_n G)/n|$$

$$\le \|G'(TC_n' TC_n - T_n' T_n)G/n\|_2$$

$$\le \|G'G\|_2 \|(TC_n' TC_n - T_n' T_n)/n\|_2.$$

Therefore, the spectral characteristics of $G'TC_n'TC_nG/n$ and those of $\Gamma^{1/2} \times X_n'X_n\Gamma^{1/2}/n = Y_2'Y_2/n$ are also asymptotically the same, by the truncation lemma, since a.s. $\Gamma^{1/2}X_n'X_n\Gamma^{1/2}/n = \Gamma^{1/2}T_n'T_n\Gamma^{1/2}/n$.

2.3.3. *Controlling the diagonal in operator norm.* Now that we have seen that a.s. we can replace $X_n$ by $TC_n$ without incurring any loss in operator norm, we will essentially focus in the analysis on the matrices obtained by replacing $X_n$ by $TC_n$ since the results will be a.s. the same. Recall that the entries of $TC_n$ are bounded by $C\sqrt{n}/(\log n)^{(1+\varepsilon)/2}$.

We turn our attention to showing that the diagonal of $G'X'XG/n$ is close to 1. We remind the reader that $G = \Sigma^{1/2}(\text{diag}(\Sigma_p))^{-1/2}$ and assume, without loss of generality, that $C \leq 2$. As a matter of fact, since $\mathbf{E}(T_n(i,i)) \to \mathbf{E}(X_n(i,i)) = 0$, $|TC_n(i,j)| \leq |X_n(i,j)| + 1$ for $n$ large enough and since $|X_n(i,j)| \leq \sqrt{n}/(\log n)^{(1+\varepsilon)/2}$, $|TC_n(i,j)| \leq 2\sqrt{n}/(\log n)^{(1+\varepsilon)/2}$ for $n$ large enough. This shows that we can take $C \leq 2$ without loss of generality.

LEMMA 4. *Let us focus on $S_p = \frac{1}{n}Y_2'Y_2 = \frac{1}{n}G'X'XG$, a quantity often studied in random matrix theory.*

*When $p \asymp n$, we have*

$$\max_{i=1,\ldots,p} |\sqrt{S_p(i,i)} - 1| \to 0 \qquad a.s.$$

PROOF. First, let $W_p = G'TC_n'TC_nG/n$. We note that, according to results in the previous subsubsection,

$$|S_p(i,i) - W_p(i,i)| = |e_i'(S_p - W_p)e_i| \leq \|S_p - W_p\|_2 \to 0 \qquad \text{a.s.}$$

Hence, the result will be shown if we can show it for $W_p(i,i)$.

We let $v_i$ denote the $i$th column of $G$. Letting $M = X'X/n$, we note that

$$S_p(i,i) = v_i'Mv_i = \|Xv_i/\sqrt{n}\|_2^2.$$

Now, consider the function $f_i$ from $\mathbb{R}^{np}$ to $\mathbb{R}$ defined by turning the vector $x$ into the matrix $X$, by first filling the rows of $X$, and then computing the Euclidean norm of $Xv_i$. In other words,

$$f_i(x) = \|Xv_i\|_2.$$

This function is clearly convex and 1-Lipschitz with respect to the Euclidean norm. As a matter of fact, for $\theta \in [0,1]$ and $x, z \in \mathbb{R}^{np}$,

$$f_i(\theta x + (1-\theta)z) = \|(\theta X + (1-\theta)Z)v_i\|_2 \leq \|\theta Xv_i\|_2 + \|(1-\theta)Zv_i\|_2$$
$$= \theta f_i(x) + (1-\theta)f_i(z).$$

Similarly,

$$|f_i(x) - f_i(z)| = |\|Xv_i\|_2 - \|Zv_i\|_2| \leq \|(X-Z)v_i\|_2 \leq \|X-Z\|_F\|v_i\|_2$$
$$= \|x-z\|_2,$$

using the Cauchy–Schwarz inequality and the fact that $\|v_i\|_2^2 = (G'G)(i,i) = \Gamma(i,i) = 1$. The same is true when $X$ is replaced by $TC_n$.

Because the $TC_n(i,j)$ are independent and bounded, we can apply recent results concerning concentration of measure of convex Lipschitz functions. In particular, from Corollary 4.10 in [33] (a consequence of Talagrand's inequality—see [46] and Theorem 4.6 in [33]), we see that for any $r > 0$, we have, if $m_{f_i}$ is a median of $f_i(TC_n)$,

$$P(|f_i(TC_n) - m_{f_i}| \geq r) \leq 4\exp(-r^2/(16C^2n/(\log n)^{(1+\varepsilon)})).$$

In particular, since $\sqrt{W_p(i,i)} = f(TC_n)/\sqrt{n}$, letting $m_{i,i}$ denote a median of $f_i(TC_n)/\sqrt{n}$, we see that

$$P(|\sqrt{W_p(i,i)} - m_{i,i}| \geq r) \leq 4\exp(-r^2(\log n)^{(1+\varepsilon)}/16C^2).$$

Finally,

$$P\Big(\Big[\max_i|\sqrt{W_p(i,i)} - m_{i,i}|\Big] \geq r\Big) \leq 4p\exp(-r^2(\log n)^{(1+\varepsilon)}/(16C^2)),$$

so, since $p \asymp n$, using the first Borel–Cantelli lemma, we see that

$$\max_{1\leq i\leq p}|\sqrt{W_p(i,i)} - m_{i,i}| \to 0 \qquad \text{a.s.}$$

All we have to do now is to show that the $m_{i,i}$ are all close to 1. We let $v_n = \text{var}(TC_n(i,j))$. Note that $v_n$ is independent of $i,j$ and that $v_n \to 1$ as $n \to \infty$.

Since we have Gaussian concentration, using Proposition 1.9 in [33], we have

$$|\mathbf{E}(\sqrt{W_p(i,i)}) - m_{i,i}| \leq 8C\sqrt{\pi}(\log n)^{-(1+\varepsilon)/2}$$

and since $\mathbf{E}(W_p(i,i)) = \|v_i\|_2^2 v_n = \Gamma(i,i)v_n = v_n$, the variance inequality in the same proposition gives

$$0 \leq v_n - \mathbf{E}(\sqrt{W_p(i,i)})^2 \leq \frac{64C^2}{(\log n)^{(1+\varepsilon)}}.$$

Consequently,

$$-\frac{8C\sqrt{\pi}}{(\log n)^{(1+\varepsilon)/2}} + \sqrt{v_n - \frac{64C^2}{(\log n)^{(1+\varepsilon)}}} \leq m_{i,i} \leq v_n + \frac{8C\sqrt{\pi}}{(\log n)^{(1+\varepsilon)/2}}.$$

Therefore, $\max_i |m_{i,i} - 1| = O(\max(|1 - v_n|, (\log n)^{-(1+\varepsilon)/2}))$ and we have

$$\max_{1 \le i \le p} |\sqrt{W_p(i,i)} - 1| \to 0 \qquad \text{a.s.}$$

We can therefore conclude that we also have

$$\max_{1 \le i \le p} |\sqrt{S_p(i,i)} - 1| \to 0 \qquad \text{a.s.} \qquad \qquad \square$$

We now turn to the more interesting situation of a covariance matrix.

LEMMA 5 (Covariance matrix).    *We now focus on the matrix*

$$S_p = \frac{1}{n-1}(Y_2 - \bar{Y}_2)'(Y_2 - \bar{Y}_2).$$

*For this matrix, we also have*

$$\max_{1 \le i \le p} |\sqrt{S_p(i,i)} - 1| \to 0 \qquad \text{a.s.}$$

PROOF.    As before, we let $W_p$ denote the equivalent of $S_p$ computed by replacing $X$ by $TC_n$.

Note that $Y_2 - \bar{Y}_2 = (\mathrm{Id}_n - \frac{1}{n}\mathbf{1}\mathbf{1}')Y_2 = (\mathrm{Id}_n - \frac{1}{n}\mathbf{1}\mathbf{1}')XG$. Now, $S_p(i,i) = v_i' \times X'(\mathrm{Id}_n - \frac{1}{n}\mathbf{1}\mathbf{1}')Xv_i/(n-1)$, so the same strategy as above can be employed, with $f$ now defined as

$$f(x) = f(X) = \left\| \left(\mathrm{Id}_n - \frac{1}{n}\mathbf{1}\mathbf{1}'\right)Xv_i \right\|_2.$$

This function is again a convex 1-Lipschitz function of $x$. Convexity is a simple consequence of the fact that norms are convex; the Lipschitz coefficient is equal to $\|v_i\|_2 \| \mathrm{Id}_n - \frac{1}{n}\mathbf{1}\mathbf{1}'\|_2$. The eigenvalues of the matrix $\mathrm{Id}_n - \frac{1}{n}\mathbf{1}\mathbf{1}'$ are $(n-1)$ ones and one zero. Its operator norm is therefore 1. We therefore have Gaussian concentration when replacing $X$ by $TC_n$. Also, we have the same bounds as before on $\max_i |S_p(i,i) - W_p(i,i)|$. All we need to check to conclude the proof is that $\mathbf{E}(W_p(i,i)) \to 1$. By renormalizing by $1/\sqrt{n-1}$, we ensure that $\mathbf{E}(W_p(i,i)) = v_n$ and so, as before, the proof is complete.    $\square$

2.3.4.  *A remark on* $\|(X - \bar{X})'(X - \bar{X})/n - 1\|_2$.    We now turn to providing a justification for Corollary 1. This amounts to understanding the behavior of the largest eigenvalue of $(Y - \bar{Y})'(Y - \bar{Y})/n - 1$, which differs slightly from what is usually investigated in the literature, namely $\widetilde{S}_p = Y'Y/n$, if, say, $Y$ is assumed to have mean zero entries.

Since in (statistical) practice $S_p = (Y - \bar{Y})'(Y - \bar{Y})/(n-1)$ is almost always used, it is of interest to know what happens for this matrix in terms of largest eigenvalue. Note that $\|S_p - \widetilde{S}_p\|_2$ does not go to zero in general, so a coarse bound of the type $|\lambda_1(S_p) - \lambda_1(\widetilde{S_p})| \leq \|\|S_p - \widetilde{S_p}\|\|_2$ is not enough to determine the behavior of $\lambda_1(S_p)$ from that of $\lambda_1(\widetilde{S_p})$.

However, letting $H_n = \mathrm{Id}_n - \frac{1}{n}11'$, we see that

$$Y - \bar{Y} = H_n Y.$$

Therefore, since $\sigma_1$, the largest singular value, is a matrix norm, we have

$$\sigma_1(Y - \bar{Y})/\sqrt{n} \leq \sigma_1(H_n)\sigma_1(Y/\sqrt{n}) = \sigma_1(Y/\sqrt{n})$$

since $H_n$ is a symmetric matrix with $(n-1)$ eigenvalues equal to 1 and one eigenvalue equal to 0.

Now, because $Y'Y/n$ and $(Y - \bar{Y})'(Y - \bar{Y})/n - 1$ have, asymptotically, the same spectral distribution, letting $l_1$ denote the right endpoint of the support of this limiting distributions (if it exists), we conclude that

$$\liminf \sigma_1^2((Y - \bar{Y})/\sqrt{n-1}) \geq l_1.$$

Hence, when $\|\|Y'Y/n\|\|_2 \to l_1$, we also have

$$\|\|(Y - \bar{Y})'(Y - \bar{Y})/(n-1)\|\|_2 \to l_1.$$

This justifies the assertion made in Corollary 1 and, more generally, the fact that when the norm of a sample covariance matrix which is not re-centered (whose entries have mean 0) converges to the right endpoint of the support of its limiting spectral distribution, so does the norm of the centered sample covariance matrix.

Finally, we note that when dealing with $S_p$, the mean of the entries of $Y$ does not matter, so we can assume without loss of generality that it is 0.

2.3.5. *Proof of Theorem 1.* We now put all of the elements together and give the proof of Theorem 1.

PROOF OF THEOREM 1. As noted in Section 2.2.1, $\mathrm{corr}(Y) = \mathrm{corr}(Y_2)$. So, to understand the spectral properties of $\mathrm{corr}(Y)$, it is enough to study those of $\mathrm{corr}(Y_2)$.

Let $S_p = (Y_2 - \bar{Y}_2)'(Y_2 - \bar{Y}_2)/(n-1)$ and $D_{S_p} = \mathrm{diag}(S_p)$. We have seen in Lemma 5 that

$$\|\|D_{S_p} - \mathrm{Id}_p\|\|_2 \to 0 \qquad \text{a.s.}$$

Now, using the remark made in Section 2.3.4, and [50], we have

$$\|\|(X - \bar{X})'(X - \bar{X})/(n-1)\|\|_2 \to (1 + \sqrt{p/n})^2 \qquad \text{a.s.}$$

Since $\|\|\Gamma_p\|\|_2$ is bounded, we have

$$\|\|S_p\|\|_2 \leq \|\|\Gamma_p\|\|_2 \|\|(X - \bar{X})'(X - \bar{X})/(n-1)\|\|_2.$$

Therefore, $\|\|S_p\|\|_2$ is a.s. bounded and the assumptions of Lemma 1 are verified.

Using Fact 1 and Lemma 1, we therefore have

$$\|\|S_p - \operatorname{corr}(Y_2)\|\|_2 \to 0 \qquad \text{a.s.}$$

Hence,

$$\|\|S_p - \operatorname{corr}(Y)\|\|_2 \to 0 \qquad \text{a.s.}$$

Finally, as explained in Section 2.3.1, the spectral properties of $S_p$, when they involve only the spectral distribution of $G'G$, are the same as those of $(Y_1 - \bar{Y}_1)'(Y_1 - \bar{Y}_1)/(n-1)$. □

**3. Elliptically distributed data and generalizations.** We now turn our attention to the problem of finding a Marčenko–Pastur-type system of equations to characterize the limiting spectral distribution of sample covariance matrices computed from elliptically distributed data and generalizations of these distributions. Our aim in doing so is manifold. From a statistical standpoint, one issue is to try and explain the lack of robustness in high dimensions of this estimate of scatter and to explain some of the numerical findings highlighted in [20]. From a more mathematical point of view, elliptical distributions raise a question concerning data vectors with a somewhat more complicated dependence structure than is usually investigated in random matrix theory. Their study will therefore force us to confront this difficulty and show that our tools allow a generalization of the results beyond elliptical distributions (and classical models).

Elliptical distributions are considered to be good models for financial data. We refer to [20] and to the book [37] for interesting discussions of the potential relevance of elliptical distributions to problems arising in the analysis of this type of data. Naturally, the study of corresponding covariance matrices is relevant to problems of portfolio optimization, where sample covariance matrices are used to estimate the covariance matrix between assets. This latter matrix is key in these problems since the optimal portfolio weights depend on the assets' covariance matrix in many formulations. Let us mention two other properties that make elliptical distributions appealing in the financial modeling context. First, we have the tail-dependence properties that they induce between components of data vectors, something that, in practice, is found in financial data and cannot be accounted for by, say, multivariate Gaussian data. Second, at least some of these distributions allow for a certain amount of "heavy-tailedness" in the observations. This is often mentioned as an important feature in modeling financial data. By contrast, it is sometimes advocated in the random matrix community that matrices with, say, i.i.d. heavy-tailed entries should be studied as models for those financial data and, in particular, returns or log-returns of stocks. We find that these more simple heavy-tail models suffer at least from one deep flaw: in the case of a crash, many companies or stocks suffer on the same day and a model of i.i.d. heavy-tailed entries does not account for this, whereas models based on elliptical distributions can. Besides the particulars of different models, it is also important to notice that the limiting spectra will be drastically different under the two types of models and the behavior of extreme eigenvalues is also very likely to be so. Before we return to our study, we refer the reader to [2] and [18] for thorough introductions to elliptical distributions.

We are therefore particularly interested in problems where we observe $n$ i.i.d. observations of an elliptically distributed vector $v$ in $\mathbb{R}^d$. Specifically, $v$ can be written as

$$v = \mu + \lambda \Gamma r,$$

where $\mu$ is a deterministic $d$-dimensional vector, $\lambda$ is a real-valued random variable, $r$ is uniformly distributed on the unit sphere in $\mathbb{R}^p$ (i.e., $\|r\|_2 = 1$) and $\Gamma$ is a $d \times p$ matrix. We let $S = \Gamma \Gamma'$. Here, $S$, a $d \times d$ matrix, is assumed to be deterministic, and $\lambda$ and $r$ are independent. We call the corresponding data matrix $X$, which is $n \times d$, that is, the vectors of observations are stacked horizontally in this matrix. Below, we will assume that $n/p$ and $d/p$ have finite limits.

As it turns out, the models for $r$ which we can handle allow for more complicated dependence structures than the one induced by a uniform distribution on the unit sphere in $\mathbb{R}^p$. So, for $r$, we will focus on random vectors whose distribution satisfies certain concentration of measure properties. For more details, we refer the reader to Section 3.2.

Note that when studying the limiting spectral distribution of a properly scaled version of $X'X/n$, we can, without loss of generality, assume that $\mu = 0$ and $\mathbf{E}(r) = 0$. As a matter of fact, if we define $\widetilde{X}$ to be the data matrix obtained by replacing $v_i$ by $\lambda_i \Gamma(r_i - \mathbf{E}(r_i))$, it is clear that $\widetilde{X}'\widetilde{X}$ is a finite rank perturbation of $X'X$ and hence properly scaled versions of these matrices have the same limiting spectral distributions. Also, since $(X - \bar{X})'(X - \bar{X})$ is a rank one perturbation of $X'X$, we see that, after proper scaling, it has the same limiting spectral distribution as the properly scaled $\widetilde{X}'\widetilde{X}$. In what follows, we will therefore assume that

$$v_i = \lambda_i \Gamma(r_i - \mathbf{E}(r_i)), \qquad i = 1, \ldots, n.$$

As is now classical, we will obtain our main result on the question of characterizing the limiting spectral distribution of a properly scaled version of $(X - \bar{X})'(X - \bar{X})$ (Theorem 2) by making use of Stieltjes transform arguments. If needed, we refer the reader to [22] for background on the connection between weak convergence of distributions and pointwise convergence of Stieltjes transforms.

We note that our model basically falls into the class of covariance matrices of the type $T_p^{1/2} X_{n,p}^* L_n X_{n,p} T_p^{1/2}$, where $X_{n,p}$ is a random matrix, independent of the square matrices $T_p$ ($p \times p$) and $L_n$ ($n \times n$), which can also be assumed to be random, as long as their spectral distributions converge to a limit. These matrices have been the subject of investigations already: see [47] Theorem 2.43, which refers to [9, 34] and [23] and the recent [41], which refers to [12] and to [51] for systems of equations involving Stieltjes transforms similar to the one we will derive. We note that under some distributional restrictions, methods of free probability using the S-transform (see [49]) could be used to derive a characterization of the limit.

However, in all of these papers, the entries of $X_{n,p}$ are assumed to be independent. Naturally, this is not the case in the situation we are considering and, as we make clear below, the dependence structures of the models we are considering can be quite complicated, as the example of the Gaussian copula (see below) indicates. Similarly, elliptically distributed data have a certain amount of dependence in the entries of the vector $r$ since its norm is 1. (We note that [35] allowed for dependence, too, and one of our questions was whether one could recover (and generalize) those results from a different angle than the one taken in [35].) Also, our matrix $\Gamma$ is $d \times p$ and usually only square matrices are considered. Interestingly, the result shows that the ratio $d/p$ plays a nontrivial role in the limiting spectral measure. One of our aims here is to show that independence in the entries of $X_{n,p}$ is not the key element; rather, we will rely on the fact that the rows of $X_{n,p}$ are independent and that the distribution of the corresponding vectors satisfies certain concentration properties.

As our proof will make clear, using the "rank one perturbation" method originally proposed in [45] and [44], proofs of convergence of spectra of random matrices basically boil down to concentration of certain quadratic forms and concentration of Stieltjes transforms, the latter being achievable using Azuma's inequality. We discuss these two aspects in Section 3.2 and Section 3.1, respectively. We chose to separate the results of these two subsections from the main proof because we believe that the results are of interest in their own right and that their technical nature would obscure the proof of Theorem 2 if they were treated there. As far as we know, many results covered by Theorem 2 are new and cannot be achieved with other methods involving (in one way or another) moment computations.

One of our points is that the importance of concentration inequalities in this context appears not to have been realized and that they permit generalizations of random matrix results to problems that look intractable by other methods.

### 3.1. *Concentration of Stieltjes transforms.*  We present a result of independent interest, namely, the fact that the Stieltjes transform of a matrix

which is the sum of $n$ independent rank one matrices is asymptotically equivalent to a deterministic function. We have somewhat more than this: we show concentration around its mean, which also immediately gives us some lower bounds on the rate of convergence.

Naturally, the result (and its extension, Remark 1) is needed in our proof (see page 37), which is why it is included here. Another reason to highlight it is the fact that it shows that certain existing results that have been obtained with "only" convergence in probability actually hold almost surely, by simply using the Borel–Cantelli lemma and the following lemma. Finally, it answers some practical questions raised in [17], which relied on Stieltjes transforms to perform spectral estimation in connection with random matrix results.

LEMMA 6 (Concentration of Stieltjes transforms). *Suppose that $M$ is a $p \times p$ matrix such that*

$$M = \sum_{i=1}^{n} r_i r_i^*,$$

*where $r_i$ are independent random vectors in $\mathbb{R}^p$. Let*

$$m_p(z) = \frac{1}{p} \operatorname{trace}((M - z \operatorname{Id}_p)^{-1}).$$

*Then, if* $\operatorname{Im}[z] = v$,

$$P(|m_p(z) - \mathbf{E}(m_p(z))| > r) \leq 4 \exp(-r^2 p^2 v^2 / (16n)).$$

Note that the lemma makes no assumptions whatsoever about the structure of the vectors $\{r_i\}_{i=1}^{n}$ other than the fact that they are independent.

PROOF OF LEMMA 6. We define $M_k = M - r_k r_k^*$. We let $\mathcal{F}_i$ denote the filtration generated by $\{r_l\}_{l=1}^{i}$. The first classical step (see [3], page 649) is to write the random variable of interest as sum of martingale differences:

$$m_p(z) - \mathbf{E}(m_p(z)) = \sum_{k=1}^{n} \mathbf{E}(m_p(z)|\mathcal{F}_k) - \mathbf{E}(m_p(z)|\mathcal{F}_{k-1}).$$

We now note that $\mathbf{E}(\operatorname{trace}((M_k - z \operatorname{Id}_p)^{-1})|\mathcal{F}_k) = \mathbf{E}(\operatorname{trace}((M_k - z \operatorname{Id}_p)^{-1})|\mathcal{F}_{k-1})$. So,

$$|\mathbf{E}(m_p(z)|\mathcal{F}_k) - \mathbf{E}(m_p(z)|\mathcal{F}_{k-1})|$$

$$= \left| \mathbf{E}(m_p(z)|\mathcal{F}_k) - \mathbf{E}\left(\frac{1}{p} \operatorname{trace}((M_k - z \operatorname{Id}_p)^{-1}) \middle| \mathcal{F}_k\right) \right.$$

$$\left. + \mathbf{E}\left(\frac{1}{p} \operatorname{trace}((M_k - z \operatorname{Id}_p)^{-1}) \middle| \mathcal{F}_{k-1}\right) - \mathbf{E}(m_p(z)|\mathcal{F}_{k-1}) \right|$$

$$\leq \left| \mathbf{E}\left( m_p(z) - \frac{1}{p}\operatorname{trace}((M_k - z\operatorname{Id}_p)^{-1})\Big|\mathcal{F}_k \right) \right|$$

$$+ \left| \mathbf{E}\left( m_p(z) - \frac{1}{p}\operatorname{trace}((M_k - z\operatorname{Id}_p)^{-1})\Big|\mathcal{F}_{k-1} \right) \right|$$

$$\leq \frac{2}{pv},$$

the last inequality following from [45], Lemma 2.6. So, $m_p(z) - \mathbf{E}(m_p(z))$ is a sum of bounded martingale differences. Note that the same would be true for its real and imaginary parts. For both of them, we can apply Azuma's inequality (see [33], Lemma 4.1) to get that

$$P(|\operatorname{Re}[m_p(z) - \mathbf{E}(m_p(z))]| > r) \leq 2\exp(-r^2 p^2 v^2/(8n)),$$

and similarly for its imaginary part. We therefore conclude that

$$P(|m_p(z) - \mathbf{E}(m_p(z))| > r) \leq P(|\operatorname{Re}[m_p(z) - \mathbf{E}(m_p(z))]| > r/\sqrt{2})$$

$$+ P(|\operatorname{Im}[m_p(z) - \mathbf{E}(m_p(z))]| > r/\sqrt{2})$$

$$\leq 4\exp(-r^2 p^2 v^2/(16n)). \qquad \square$$

We have the following, immediate, corollary.

COROLLARY 2. *Suppose that we consider the following sequence of random matrices: for each $p$, select $n$ independent $p$-dimensional vectors. Let $M = \sum_{i=1}^n r_i r_i^*$. Assume that $p/n$ remains bounded away from 0. Then,*

$$\forall z \in \mathbb{C}^+ \qquad m_p(z) - \mathbf{E}(m_p(z)) \to 0 \qquad a.s.$$

*and also*

$$\forall z \in \mathbb{C}^+ \qquad \frac{\sqrt{p}}{(\log p)^{(1+\alpha)/2}}|m_p(z) - \mathbf{E}(m_p(z))| \to 0 \qquad a.s., \text{ for } \alpha > 0.$$

*In other words, $m_p(z)$ is asymptotically deterministic.*

PROOF. The proof is an immediate consequence of the first Borel–Cantelli lemma. $\square$

REMARK 1. We note that if $\Sigma$ is a matrix independent of the $r_i$, similar results would apply to

$$\frac{1}{p}\operatorname{trace}((M - z\operatorname{Id}_p)^{-1}\Sigma^l)$$

after we replace $v$ by $v/\|\Sigma\|_2^l$. In particular, if $\|\Sigma\|_2 \leq C(\log p)^m$ for some $m$, we have

$$\frac{1}{p}\operatorname{trace}((M - z\operatorname{Id}_p)^{-1}\Sigma^l) - \mathbf{E}\left(\frac{1}{p}\operatorname{trace}((M - z\operatorname{Id}_p)^{-1}\Sigma^l)\right) \to 0 \qquad a.s.$$

However, the rate in the second part of the previous corollary needs to be adjusted.

REMARK 2. We note that the rate given by Azuma's inequality does not match the rate that appears in results concerning fluctuation behavior of linear spectral statistics, which is $n$ and not $\sqrt{n}$. Of course, our result encompasses many situations that are not covered by the currently available results on linear spectral statistics, which might help to explain this discrepancy. The "correct" rate can be recovered using ideas similar to the ones discussed in [26] and [33], Chapter 8, Section 5. As a matter of fact, if we consider the Stieltjes transform of the measure that puts mass $1/p$ at each of the singular values of $M = X^*X/n$, it is an easy exercise to see that this function (of $X$) is $1/(\sqrt{np}v^2)$-Lipschitz with respect to the Euclidean (or Frobenius) norm. Hence, if the $np$-dimensional vector made up of the entries of $X$ has a distribution that satisfies a dimension-free concentration property with respect to the Euclidean norm, we find that the fluctuations of the Stieltjes transform at $z$ are of order $\sqrt{np}$, which corresponds to the "correct" rate found in the analysis of these models. (Note, however, that results have been shown beyond the case of distributions with dimension-free concentration.)

The conclusion of this discussion is that since the spectral distribution of random matrices is characterized by their Stieltjes transforms, it is not surprising that they are asymptotically nonrandom for a very wide class of data matrices of covariance type. We now turn to the examination of another type of concentration which we will crucially need in our proof, namely the concentration of certain quadratic forms.

3.2. *Concentration of quadratic forms.* The key property we will rely on in the proof of our main theorem is a concentration property of quadratic forms. This property is summarized in Corollary 4 and we now give important sufficient conditions to reach it.

LEMMA 7 (Case of Gaussian concentration). *Suppose that the random vector $r \in \mathbb{R}^p$ has the property that for any convex 1-Lipschitz (with respect to the Euclidean norm) function $F$ from $\mathbb{R}^p$ to $\mathbb{R}$, we have, if $m_F$ denotes a median of $F(r)$,*

$$P(|F(r) - m_F| > t) \le C \exp(-c(p)t^2),$$

*where $C$ and $c(p)$ are independent of $F$ and $C$ is independent of $p$. We allow $c(p)$ to be a constant or to go to zero with $p$ like $p^{-\alpha}$, $0 \le \alpha < 1$. Suppose, further, that $\mathbf{E}(r) = 0$, $\mathbf{E}(rr^*) = \Sigma$, with $\|\Sigma\|_2 \le \log(p)$.*

*If $M$ is a complex deterministic matrix such that $|\!|\!|M|\!|\!|_2 \leq \xi$, where $\xi$ is independent of $p$, then*

$$\frac{1}{p}r'Mr \text{ is strongly concentrated around its mean, } \frac{1}{p}\operatorname{trace}(M\Sigma).$$

*In particular, if, for $\varepsilon > 0$, $t_p(\varepsilon) = \log(p)^{1+\varepsilon}/\sqrt{pc(p)}$, then*

$$\log\left\{P\left(\left|\frac{1}{p}r'Mr - \frac{1}{p}\operatorname{trace}(M\Sigma)\right| > t_p(\varepsilon)\right)\right\} \asymp -(\log p)^{1+2\varepsilon}.$$

*If $\mathbf{E}(r) \neq 0$, then the same results are true when one replaces $r$ by $r - \mathbf{E}(r)$ everywhere and $\Sigma$ is the covariance of $r$.*

*Finally, if $\xi$ is allowed to vary with $p$, then the same results hold when one replaces $t_p(\varepsilon)$ by $t_p(\varepsilon)\xi$ or, equivalently, divides $M$ by $\xi$.*

PROOF. In what follows, $K$ denotes a generic constant that may change from occurrence to occurrence, but which is independent of $p$. First, it is clear that we can rewrite $M$ as $M = RM + iIM$, where $RM$ and $IM$ are real matrices. Further, the spectral norm of those matrices is less than $\xi$ [of course, $RM = (M + M_1)/2$, where $M_1$ is the (entrywise) complex conjugate of $M$].

Now, strong concentration for $r'RMr/p$ and $r'IMr/p$ will imply strong concentration for the sum of those two terms. We note that since $r'RMr$ is real, $r'RMr = (r'RMr)'$ and

$$r'RMr = r'\left(\frac{RM + RM'}{2}\right)r.$$

Hence, instead of working on $RM$, we can work on its symmetrized version.

Let us now decompose $(RM + RM')/2$ into $RM_+ + RM_-$, where $RM_+$ is positive semidefinite and $-RM_-$ is positive definite [or 0 if $(RM + RM')/2$ itself is positive semidefinite]. This is possible because $(RM + RM')/2$ is real symmetric and we carry out this decomposition by simply following its spectral decomposition. Note that both matrices have spectral norm less than $\xi$. Now, the map $\phi : r \to \sqrt{r'RM_+r/p}$ is $\sqrt{\xi/p}$-Lipschitz (with respect to the Euclidean norm) and convex, which is easily seen after one notes that $\sqrt{r'RM_+r/p} = \|RM_+^{1/2}r/\sqrt{p}\|_2$. This guarantees, by our assumption, that

$$P(|\sqrt{r'RM_+r/p} - m_\phi| > t) \leq C\exp(-pc(p)t^2/\xi),$$

where $m_\phi$ is a median of $\phi(r)$.

Now, using Proposition 1.9 in [33], letting $\mu_\phi$ denote the mean of $\phi(r)$ and observing that $\operatorname{var}(\phi(r)) = \mathbf{E}(r'RM_+r/p) - \mu_\phi^2$, we have

$$|\mu_\phi - m_\phi| \leq \frac{C\sqrt{\pi}}{2}\sqrt{\frac{1}{pc(p)}} \quad \text{and}$$

$$0 \leq \mathbf{E}(r'RM_+r/p) - \mu_\phi^2 \leq \frac{C}{pc(p)}.$$

We hence deduce, using the fact that $\sqrt{a+b} \leq \sqrt{a} + \sqrt{b}$ for nonnegative reals, that

$$|\sqrt{\mathbf{E}(r'RM_+r/p)} - m_\phi| \leq \kappa_p = \frac{C\sqrt{\pi}/2 + \sqrt{C}}{\sqrt{pc(p)}}.$$

We also have, trivially,

$$\mathbf{E}(r'RM_+r/p) = \frac{\mathrm{trace}(RM_+\Sigma)}{p}$$

since $\mathbf{E}(r) = 0$. Therefore, we have that for $u > 0$,

$$P(|\sqrt{r'RM_+r/p} - \sqrt{\mathrm{trace}(RM_+\Sigma/p)}| > u + \kappa_p) \leq C \exp(-pc(p)u^2/\xi).$$

On the other hand, if $0 < t < \kappa_p$, then we have

$$P(|\sqrt{r'RM_+r/p} - \sqrt{\mathrm{trace}(RM_+\Sigma/p)}| > t)$$
$$\leq 1 \leq \exp(pc(p)\kappa_p^2/\xi)\exp(-pc(p)(t-\kappa_p)^2/\xi)$$
$$\leq \max(C, \exp(pc(p)\kappa_p^2/\xi))\exp(-pc(p)(t-\kappa_p)^2/\xi).$$

Since $pc(p)\kappa_p^2 = (C\sqrt{\pi}/2 + \sqrt{C})^2$, we conclude that for any $t > 0$,

$$P(|\sqrt{r'RM_+r/p} - \sqrt{\mathrm{trace}(RM_+\Sigma/p)}| > t) \leq K \exp(-pc(p)(t-\kappa_p)^2/\xi),$$

where $K$ depends on $C$ and $\xi$, but not on $p$. Note that $\kappa_p \to 0$ since $pc(p) \to \infty$ and $C$ does not depend on $p$.

We now turn to finding a deviation inequality for the quadratic form of interest. Let us define

$$\zeta_p = \mathrm{trace}(RM_+\Sigma/p),$$
$$A = \{|r'RM_+r/p - \zeta_p| > t\},$$
$$B = \{\sqrt{r'RM_+r/p} \leq \sqrt{\zeta_p} + 1\}.$$

Our aim is to show that the probability of $A$ is "exponentially small" in $p$. Of course, we have $P(A) \leq P(A \cap B) + P(B^c)$. We note that $P(B^c)$ is "exponentially small" in $p$ since

$$P(B^c) = P(\sqrt{r'RM_+r/p} - \sqrt{\mathrm{trace}(RM_+\Sigma/p)} \geq 1) \leq K \exp(pc(p)(1-\kappa_p)^2/\xi)$$

and $pc(p) \to \infty$, at least as fast as $p^{1-\alpha}$, with $\alpha > 0$. Now, note that

$$A \cap B \subseteq D = \left\{|\sqrt{r'RM_+r/p} - \sqrt{\zeta_p}| > \frac{t}{2\sqrt{\zeta_p} + 1}\right\}.$$

To see this, note simply that for positive reals, $|x - y| = |\sqrt{x} - \sqrt{y}|(\sqrt{x} + \sqrt{y})$. Finally, because of our bounds on the norm of $\Sigma$ and the fact that $\|RM_+\|_2 \leq \xi$, we see that $\text{trace}(RM_+\Sigma/p) = \zeta_p \leq \log(p)\xi$. Hence, $P(D) \leq K \exp(-pc(p)(t/(2\sqrt{\zeta_p} + 1) - \kappa_p)^2/\xi)$ for some $K$ independent of $p$ and we therefore have

$$P(A) \leq K[\exp(-pc(p)(t/(2\sqrt{\zeta_p} + 1) - \kappa_p)^2/\xi) + \exp(-pc(p)(1 - \kappa_p)^2/\xi)]$$

$$= g_p(t).$$

Similarly, we can obtain the same type of bounds for $\sqrt{-r'RM_-r/p}$. From those, we conclude that

$$P(|r'RMr/p - \text{trace}(RM\Sigma_p)/p| > t) \leq 2g_p(t/2).$$

Finally,

$$P(|r'Mr/p - \text{trace}(M\Sigma_p)/p| > t) \leq 4g_p(t/2\sqrt{2}).$$

From the way $g_p(t)$ behaves, we obtain the result concerning strong concentration. Now, studying the asymptotics of $g_p(t_p(\varepsilon))$ for large $p$ gives the statement concerning the log probability in the lemma.

The last statement in the lemma, concerning the replacement of $r$ by $r - \mathbf{E}(r)$ if $\mathbf{E}(r) \neq 0$, follows simply from the fact that, for any given $\mu$, the map $\phi(r) = \sqrt{(r - \mu)'\Sigma(r - \mu)/p}$ is convex and $\sqrt{\xi/p}$-Lipschitz since the composition of a convex mapping and an affine one is convex. (See, e.g., Section 3.2.2 in [10].)  □

Motivated by examples we will see below, we also note that we have the following corollary, which is applicable when concentration is not limited to convex 1-Lipschitz functions, but holds for any 1-Lipschitz function.

COROLLARY 3 (Gaussian concentration of nonconvex functionals). *Suppose that the random vector $v \in \mathbb{R}^p$ has the property that for any 1-Lipschitz (with respect to the Euclidean norm) function $F$ from $\mathbb{R}^p$ to $\mathbb{R}$, we have, if $m_F$ denotes a median of $F(v)$,*

$$P(|F(v) - m_F| > t) \leq C \exp(-c(p)t^2),$$

*where $C$ and $c(p)$ are independent of $F$ and $C$ is independent of $p$. We allow $c(p)$ to be a constant or to go to zero with $p$ like $p^{-\alpha}$, $0 \leq \alpha < 1$.*

*Consider $r = \Phi(v)$, where $\Phi$ is a 1-Lipschitz map from $\mathbb{R}^p$ to $\mathbb{R}^d$, also with respect to the Euclidean norm. Suppose, as above, that $\mathbf{E}(r) = 0$, $\mathbf{E}(rr^*) = \Sigma$, with $\|\Sigma\|_2 \leq \log(p)$.*

*If $M$ is a complex deterministic matrix such that $\|M\|_2 \leq \xi$, where $\xi$ is independent of $p$, then*

$$\frac{1}{p}r'Mr \text{ is strongly concentrated around its mean, } \frac{1}{p}\text{trace}(M\Sigma).$$

*In particular, if for $\varepsilon > 0$, $t_p(\varepsilon) = \log(p)^{1+\varepsilon}/\sqrt{pc(p)}$, then*

$$\log\left\{ P\left(\left|\frac{1}{p}r'Mr - \frac{1}{p}\operatorname{trace}(M\Sigma)\right| > t_p(\varepsilon)\right)\right\} \asymp -(\log p)^{1+2\varepsilon}.$$

*If $\mathbf{E}(r) \neq 0$, then the same is true when one replaces $r$ by $(r - \mathbf{E}(r))$ everywhere and $\Sigma$ is the covariance of $r$. Finally, if $\xi$ is allowed to vary with $p$, then the same results hold when one replaces $t_p(\varepsilon)$ by $t_p(\varepsilon)\xi$ or, equivalently, divides $M$ by $\xi$.*

PROOF.    The proof follows easily from the arguments developed for Lemma 7 after we note that the map $\phi \colon R^p \to R$ with $\phi(v) = \sqrt{r'RM_+r/p} = \sqrt{\Phi(v)'RM_+\Phi(v)/p}$ is $\sqrt{\xi/p}$-Lipschitz with respect to the Euclidean norm. The concentration properties of $v$ can then be invoked and the proof follows along similar lines as above.    □

For applications, it is important to extend the results beyond Gaussian concentration. We therefore state the following lemma.

LEMMA 8 (Beyond Gaussian concentration).    *Suppose that the random vector $r \in \mathbb{R}^p$ has the property that for any convex 1-Lipschitz (with respect to the Euclidean norm) function $F$ from $\mathbb{R}^p$ to $\mathbb{R}$, we have, for $b > 0$ independent of $p$ and $m_F$ a median of $F$,*

$$P(|F(r) - m_F| > t) \leq C \exp(-c(p)t^b),$$

*where $C$ and $c(p)$ are independent of $F$ and $C$ is independent of $p$. We allow $c(p)$ to be a constant or to go to zero with $p$ like $p^{-\alpha}$, $0 \leq \alpha < b/2$. Suppose, further, that $\mathbf{E}(r) = 0$, $\mathbf{E}(rr^*) = \Sigma$, with $\|\Sigma\|_2 \leq \log(p)$.*

*If $M$ is a complex deterministic matrix such that $\|M\|_2 \leq \xi$, where $\xi$ is independent of $p$, then*

$$\frac{1}{p}r'Mr \text{ is strongly concentrated around its mean, } \frac{1}{p}\operatorname{trace}(M\Sigma).$$

*In particular, if, for $\varepsilon > 0$, $t_p(\varepsilon) = \log(p)^{1/2+1/b+\varepsilon}/\sqrt{pc^{2/b}(p)}$, then*

$$\log\left\{ P\left(\left|\frac{1}{p}r'Mr - \frac{1}{p}\operatorname{trace}(M\Sigma)\right| > t_p(\varepsilon)\right)\right\} \asymp -(\log p)^{1+b\varepsilon}.$$

*If $\mathbf{E}(r) \neq 0$, then the same is true when one replaces $r$ by $(r - \mathbf{E}(r))$ everywhere and $\Sigma$ is the covariance of $r$.*

*Finally, if $\xi$ is allowed to vary with $p$, then the same results hold when one replaces $t_p(\varepsilon)$ by $t_p(\varepsilon)\xi$ or, equivalently, divides $M$ by $\xi$.*

Proof. We only give a sketch of the proof. The ideas are exactly the same as those above. However, when studying the concentration of $\sqrt{r'RM_+r/p}$, the exponent of the exponential is, to leading order, $p^{b/2}c(p)(t - \kappa_p)^b$. We note that $\kappa_p$ will be somewhat different in its form than it was in the Gaussian concentration case. This comes from the fact, following the analysis in Proposition 1.9 of [33], that the inequalities we now have, if $\mu_F$ denotes the mean of $F$, are

$$|\mu_F - m_F| \leq \frac{C}{bc^{1/b}}\Gamma\left(\frac{1}{b}\right) \quad \text{and} \quad \text{var}(F) \leq \frac{2C}{bc^{2/b}}\Gamma\left(\frac{2}{b}\right),$$

where $\Gamma$ denotes the Gamma function. With this adjustment, the previous proof proves the present lemma. □

A corollary similar to Corollary 3 holds for the variant of Lemma 8 where concentration is not limited to convex 1-Lipschitz functionals, but is valid for any 1-Lipschitz (with respect to the Euclidean norm) functional.

*Examples of distributions for which the previous results apply.*

- Gaussian random variables with $|||\Sigma|||_2 \leq \log(p)$. Lemma 7 and Corollary 3 apply, according to [33], Theorem 2.7, with $c(p) = 1/|||\Sigma|||_2$.
- Vectors of the type $\sqrt{p}r$, where $r$ is uniformly distributed on the unit ($\ell_2$-) sphere in dimension $p$. Theorem 2.3 in [33] shows that Lemma 7 (and Corollary 3) applies, with $c(p) = (1 - 1/p)/2$, after noting that a 1-Lipschitz function with respect to the Euclidean norm is also 1-Lipschitz with respect to the geodesic distance on the sphere. As we will see below, this will allow us to treat the case of elliptically distributed data.
- Vectors $\Gamma\sqrt{p}r$, with $r$ uniformly distributed on the unit ($\ell_2$-) sphere in $\mathbb{R}^p$ and with $\Gamma\Gamma' = \Sigma$ having the characteristics explained in Lemma 7.
- Vectors of the type $p^{1/b}r$, $1 \leq b \leq 2$, where $r$ is uniformly distributed in the unit $\ell^b$ ball or sphere in $\mathbb{R}^p$. (See [33], Theorem 4.21, which refers to [43] as the source of the theorem.) Lemma 8 applies to them, with $c(p)$ depending only on $b$.
- Vectors with log-concave density of the type $e^{-U(x)}$, with the Hessian of $U$ satisfying, for all $x$, $\text{Hess}(U) \geq c\,\text{Id}_p$, where $c > 0$ has the characteristics of $c(p)$ in Lemma 7; see [33], Theorem 2.7. Here, we also need $|||\Sigma|||_2$ to satisfy the assumptions of Lemma 7. Corollary 3 also applies here.
- Vectors $(r)$ distributed according to a (centered) Gaussian copula, with corresponding correlation matrix $\Sigma$ such that $|||\Sigma|||_2$ is bounded. Here, we can apply Corollary 3 since if $\tilde{r}$ has a Gaussian copula distribution, then its $i$th entry satisfies $\tilde{r}_i = \Phi(v_i)$, where $v$ is multivariate normal with covariance matrix $\Sigma$, $\Sigma$ being a correlation matrix, that is, its diagonal is 1. Here, $\Phi$ is the cumulative distribution function of a standard normal distribution, which is trivially Lipschitz. Now, taking $r = \tilde{r} - 1/2$ gives a centered Gaussian copula. The fact that the covariance matrix of $r$ then has bounded operator norm requires a little work and is shown in the Appendix.

- Vectors with i.i.d. entries bounded by $1/\sqrt{c(p)}$. See Corollary 4.10 in [33] for the concentration part, which shows that Lemma 7 applies. We crucially need the fact that the concentration of measure result is valid "only" for convex 1-Lipschitz functions. As we will explain below, in our main theorem, this result will enable us to work with random variables with bounded second moment since using an argument similar to those in [45], those random variables can be truncated at $\log(p)$ without (a.s.) changing the limiting spectral distribution.

It also appears possible to use this method to treat more "exotic" examples involving vectors sampled uniformly from certain Riemannian submanifolds of $\mathbb{R}^p$, a question which is sometimes of interest in multivariate statistics and computer science. We refer to [33], Theorems 2.4 and 3.1, for the concentration aspects of these questions.

We also have the following, important, corollary that will play a key role in the proof of Theorem 2.

COROLLARY 4. *Suppose that $\{r_i\}_{i=1}^n$ are independent random vectors whose distributions satisfy the hypotheses of Lemmas 7 or 8, or Corollary 3. Suppose that $n \asymp p$. Suppose that $M_i$ are random matrices, $M_i$ being independent of $r_i$ and such that $\|M_i\|_2 \le K$, where $K$ is nonrandom. Suppose, further, that for some matrix $\mathcal{M}$, some $K_p$ with $K_p = O(K/p)$ and all Hermitian matrices $A$,*

$$\forall i \qquad \left| \frac{1}{p} \operatorname{trace}(M_i A) - \frac{1}{p} \operatorname{trace}(\mathcal{M} A) \right| \le \|A\|_2 K_p.$$

*Then, for any $\varepsilon > 0$, if $\tilde{r}_i = r_i - \mathbf{E}(r_i)$, we have*

$$(1) \quad \frac{\sqrt{p c^{2/b}(p)}}{(\log p)^{(1/2 + 1/b + \varepsilon)} K} \max_{i=1,\dots,n} \left| \frac{1}{p} \tilde{r}_i' M_i \tilde{r}_i - \frac{1}{p} \operatorname{trace}(\mathcal{M} \Sigma) \right| \to 0 \qquad a.s.$$

PROOF. From the previous results, we have

$$P\left( \max_i \left| \frac{1}{p} \tilde{r}_i' M_i \tilde{r}_i - \frac{1}{p} \operatorname{trace}(M_i \Sigma) \right| > t \right)$$

$$\le \sum_{i=1}^n P\left( \left| \frac{1}{p} \tilde{r}_i' M_i \tilde{r}_i - \frac{1}{p} \operatorname{trace}(M_i \Sigma) \right| > t \right)$$

$$\le 4n g_p(t/2\sqrt{2}),$$

by conditioning on $M_i$ to compute each probability in the sum. Therefore, using the first Borel–Cantelli lemma, the results above and the fact that $n \asymp p$, we have

$$\frac{\sqrt{pc^{2/b}(p)}}{(\log p)^{(1/2+1/b+\varepsilon)} K} \max_i \left| \frac{1}{p} \tilde{r}_i' M_i \tilde{r}_i - \frac{1}{p} \operatorname{trace}(M_i \Sigma) \right| \to 0 \qquad \text{a.s.}$$

and because $|\frac{1}{p} \operatorname{trace}(M_i \Sigma) - \frac{1}{p} \operatorname{trace}(\mathcal{M} \Sigma)| \leq K_p \|\Sigma\|_2 \leq K_p \log(p)$, we conclude that

$$\frac{\sqrt{pc^{2/b}(p)}}{(\log p)^{(1/2+1/b+\varepsilon)} K} \max_i \left| \frac{1}{p} \tilde{r}_i' M_i \tilde{r}_i - \frac{1}{p} \operatorname{trace}(\mathcal{M} \Sigma) \right| \to 0 \qquad \text{a.s.} \qquad \square$$

We also have the following technical result that will be useful below.

COROLLARY 5. *We assume that Lemma 7, Lemma 8 or Corollary 3 applies and that* $\operatorname{trace}(\Sigma)/p$ *is bounded by $K$ independent of $p$. If $\{r_i\}_{i=1}^n$ is a triangular array of independent random variables and $n/p$ remains bounded, then the spectral distribution of $A_n = \sum_{i=1}^n r_i r_i^* / n$ is a.s. tight.*

PROOF.   We first assume that $\mathbf{E}(r_i) = 0$. Let $R_n$ denote the matrix whose $i$th row is $r_i'$. We consider the first moment of the spectral distribution of $A_n = R_n^* R_n / n$, which is equal to $M_1$, with $M_1 = 1/n \sum_{i=1}^n \operatorname{trace}(r_i r_i^* / p)$. Its mean is $\operatorname{trace}(\Sigma)/p$. As we just saw, $r_i^* r_i / p$ is strongly concentrated around $\operatorname{trace}(\Sigma)/p(= \mathbf{E}(M_1))$ and this property transfers to $M_1$ using the fact that $P(|M_1 - \mathbf{E}(M_1)| > t) \leq n P(|r_i^* r_i / p - \mathbf{E}(r_i^* r_i / p)| > t)$. Because $\operatorname{trace}(\Sigma)/p$ is assumed to be bounded, we see that $M_1$ is a.s. bounded by $K + 1$. Because it is the first moment of the spectral distribution of $A_n$, we conclude that, if we let $F^{A_n}$ denote the c.d.f. of the spectral distribution of $A_n$, we have a.s. $F^{A_n}([M, \infty)) \leq (K+1)/M$ for $n$ sufficiently large. Since the spectral distribution of $A_n$ is supported on $[0, \infty)$, we conclude that it is a.s. tight.

In the case where $r_i$ do not have mean 0, we can work with $\tilde{r}_i = r_i - \mathbf{E}(r_i)$. The resulting matrix $\tilde{R}_n$ is a perturbation of $R_n$ of rank at most 3 and therefore generates the same limiting spectral distribution as that generated by $R_n$. Therefore, the previous arguments applied to $\tilde{R}_n$ give the result for $R_n$.  $\square$

We conclude this concentration discussion by considering some practical consequences.

*Practical geometric consequences of concentration.*   If we apply the previous results to the matrix $M = \operatorname{Id}_p$, we see that our concentration results indicate that if $\mathbf{E}(r) = 0$, then $\|r\|^2/p$ is strongly concentrated around

trace$(\Sigma_p)/p$. In particular, if $n \asymp p$, we see that $\max_{i=1,\dots,n}|\,\|r\|^2/p - \mathrm{trace}(\Sigma_p)/p|$ will tend to 0 a.s. Hence, the vectors $r_i/\sqrt{p}$ appear to be located close to a sphere. By properly choosing $M$, for instance, as a block matrix with 0 on the block diagonal and $\mathrm{Id}_p$ off the diagonal, we can also show (see [15] for details) that $\max_{i \neq j}|r_i' r_j/p|$ tends to 0 a.s. and hence so does the maximal angle between two such vectors. So, perhaps surprisingly, the vectors $r_i$ appear to be almost orthogonal to one another. These remarks suggest that although the models we study are quite general, their geometric properties are somewhat peculiar. Hence, one should probably check (even just graphically) whether such geometric features are present in the data before using random matrix results.

### 3.3. *Marčenko–Patur-type system for covariance matrices computed from generalized elliptically distributed data.*
We refer the reader to the discussion introducing Section 3 for a review of literature concerning elliptically distributed data and some motivation for the theorem that follows. In what follows, we assume that we have a triangular "array" of random variables, where the $n$th line contains $n$ i.i.d. $\lambda_i$'s and $n$ i.i.d. $r_i$'s in $\mathbb{R}^p$ satisfying concentration inequalities as in Lemma 7, Lemma 8 or Corollary 3. We assume that the $r_i$'s have covariance matrix $\Sigma$ such that $\|\Sigma\|_2 \leq \log(p)$. We also have to work with a $d \times p$ matrix $\Gamma$. The data vectors we will focus on are therefore the "array" of

$$v_i = \mu + \lambda_i \Gamma r_i,$$

which we say have generalized elliptical distributions. In what follows, we allow $S = \Gamma\Gamma'$ to be random, as long as it is independent of the vectors $r_i$. For all practical purposes, however, $S$ can be considered deterministic.

We present the theorem in the form that makes it most natural for elliptically distributed data, our original motivation.

**Theorem 2.** *Let $\{\{v_i\}_{i=1}^{n}\}_{n=1}^{\infty}$ form a triangular array of independent random vectors, "generalized-elliptically" distributed, as described above. In particular, recall that they are in $\mathbb{R}^d$.*

- *Define $\theta_n = d/p$, $\rho_n = p/n$, $\xi_n = d^2/np = \theta_n^2 \rho_n$.*
- *Let $G_d$ denote the spectral distribution of $\Gamma\Gamma' = S$, $H_d$ the spectral distribution of $\Gamma\Sigma\Gamma' = T$ ($S$ and $T$ are $d \times d$) and $\nu_n$ the spectral distribution of the diagonal matrix containing the $\lambda_i$'s.*
- *Assume that $H_d$ converges weakly a.s. to a probability distribution $H \neq 0$. Assume, further, that $\int \tau \, dH_d(\tau)$ remains bounded.*
- *Assume that $G_d$ converges weakly a.s. to a probability distribution $G \neq 0$.*
- *Assume that $\nu_n$ converges weakly a.s. to a probability distribution $\nu \neq 0$.*

*Let $X$ denote the $n \times d$ data matrix whose $i$th row is $v_i$. Consider the matrix*

$$B_n = \frac{d}{p}\frac{1}{n}X'X = \frac{\theta_n}{n}\sum_{i=1}^n v_i v_i' \triangleq \sum_{i=1}^n u_i u_i'.$$

*If $\rho_n$ has a finite nonzero limit $\rho$ and $\theta_n$ has a finite nonzero limit $\theta$, then $\xi_n$ obviously has a finite nonzero limit $\xi$ and the Stieltjes transform of $B_n$, $m_n$, converges a.s. to a deterministic limit $m$ satisfying the equations*

$$m(z) = \int \frac{dH(\tau)}{\tau \int \theta\lambda^2/(1+\xi\lambda^2 w(z))\, d\nu(\lambda) - z} \qquad and$$

$$w(z) = \int \frac{\tau\, dH(\tau)}{\tau \int \theta\lambda^2/(1+\xi\lambda^2 w(z))\, d\nu(\lambda) - z}.$$

*$w$ is the unique solution of this equation mapping $\mathbb{C}^+$ into $\mathbb{C}^+$. (The intuitive meaning of $w$ is explained below. We also remind the reader that $m$ uniquely characterizes the limiting spectral distribution of $B_n$.)*

*We note, further, that we have*

$$1 + zm(z) = w(z)\int \frac{\theta\lambda^2}{1+\xi\lambda^2 w(z)}\, d\nu(\lambda).$$

*The same results hold for the scaled sample covariance matrix $d/p(X - \bar{X})'(X - \bar{X})/n$ since it is a finite-rank perturbation of $B_n$.*

*The conclusion is that the limiting spectral distribution of $B_n$ is a nonrandom probability measure and is uniquely characterized by the previous system of two equations.*

We note that our finite-rank perturbation arguments (see the introduction to Section 3) allow $\mu$ and $\mathbf{E}(r)$ to be arbitrary. However, in the proof, we can, and will, assume (without loss of generality) that $\mu = 0$ and $\mathbf{E}(r) = 0$.

In the proof, we do not actually need the $\lambda_i$'s to be independent of each other. We only need them to be independent of the $r$'s and their empirical distribution to converge a.s. to a deterministic limit, $\nu$. In the case of i.i.d. $\lambda_i$'s, we note that $\nu_n$ has an almost sure limit $\nu$ by the Glivenko–Cantelli theorem ([48], Theorem 19.1) for triangular arrays. (A simple modification to the proof given in [48], which is not for triangular arrays, can be obtained using Hoeffding's inequality for the variables $1_{\lambda_i \le t}$, which guarantees that the result is true for triangular arrays.)

We note that, perhaps interestingly, the proof could be adapted to show that quantities of the type $\operatorname{trace}(T^k(B_n - z\operatorname{Id}_p)^{-1})/d$ satisfy the same equation as $w$, with $\tau$ raised to the power $k$ at the numerator and the same denominator involving $w$, provided the $H_d$'s have enough moments. (Note that this is the case for $m$, with $k = 0$ and $w$ which basically corresponds to $k = 1$.)

To make the theorem more concrete, we now give a few examples of distributions to which it can be applied. The concentration justifications appear in Section 3.2.

- Elliptical distributions. In this case, $r_i = \sqrt{p}\tilde{r}_i$, with $\tilde{r}_i$ uniformly distributed on the sphere, so $\Sigma = \mathrm{Id}_p$. Note, in particular, that $\lambda_i$ can have a Cauchy distribution or any heavy-tailed distribution. The theorem hence describes the limiting measure obtained when using data sampled according to the multivariate $t$ or Cauchy distributions.

- Data distributed according to a Gaussian copula, with corresponding correlation matrix, $R$, bounded in operator norm. In this case, $\lambda_i = 1$, $r_i = \Phi(\tilde{r}_i)$, where $\tilde{r}_i = \mathcal{N}(0, R)$, and $\Phi$ is the c.d.f. of the standard normal distribution. The theorem then says that the Marčenko–Pastur equation holds when $r_i$ is sampled according to this distribution. This example, in particular, appears to be out of reach of methods relying in one way or another on moment computations.

- $r_i = p^{1/b}\tilde{r}_i$, where $\tilde{r}_i$ is sampled uniformly from the unit $\ell_b$-ball or sphere, $1 \leq b \leq 2$, in $\mathbb{R}^p$. We refer the reader to [33] pages 37–38 for some of the subtleties which arise for the sphere when $1 < b < 2$.

- $r_i$ has i.i.d. entries with finite second moment. Then, using the truncation arguments in [45], we see that we can truncate $r_i$ at level $\log(p)$ without a.s. affecting the limiting spectral distribution. (The arguments in [45] are rank arguments and carry directly over to our situation.) We then have $c(p) = (\log(p))^{-2}$ when using Lemma 7. Here, the convexity assumption mentioned in Lemma 7 is necessary, as we rely crucially on Corollary 4.10 in [33] for the concentration arguments.

- We note that if $\lambda_i = 1$ and $d = p$, then we recover the Marčenko–Pastur equation. The theorem therefore provides an extension of the known range of validity of this result. (We note that our result is in that case related to [39].) In this setting, the practical geometric remark made at the end of Section 3.2 applies and, hence, one should probably perform simple graphical diagnostics on the data before relying on insights drawn from random matrix results.

The system of equations we have found is unfortunately not trivial to exploit in order to gain further understanding of the spectra of the matrices at stake; we postpone a detailed investigation of its consequences to a further project. We now turn to the proof of Theorem 2.

3.3.1. *Preliminaries.* We note that the matrix we are considering is of the form $\Gamma X' DX \Gamma'$, where $D$ is a diagonal matrix containing the $\lambda_i^2$'s, that is, $D_{i,j} = 1_{i=j}\lambda_i^2$. We let $s_i$ denote the eigenvalues of $S = \Gamma\Gamma'$.

We denote by $\|F\|$ the value $\sup_x |F(x)|$ and by $F^M$ the c.d.f. of the spectral distribution of the matrix $M$. We see, using Lemma 2.5 in [45],

that

$$\|F^{Q^*T_1Q} - F^{\tilde{Q}^*\tilde{T}_1\tilde{Q}}\| \leq \frac{1}{p}(\operatorname{rank}(T_1 - \tilde{T}_1) + 2\operatorname{rank}(\tilde{Q} - Q)).$$

In our situation, we have $Q = X\Gamma'$ and $T_1 = D$, so, using the fact that $\operatorname{rank}(AB) \leq \min(\operatorname{rank}(A), \operatorname{rank}(B))$, we conclude that

$$\|F^{Q^*T_1Q} - F^{\tilde{Q}^*\tilde{T}_1\tilde{Q}}\| \leq \frac{1}{p}(\operatorname{rank}(D - \tilde{D}) + 2\operatorname{rank}(\tilde{\Gamma}' - \Gamma')).$$

Let us now choose for $\tilde{D}$ the diagonal matrix with entries $\lambda_i^2 1_{\lambda_i^2 \leq \alpha_p}$, which we abbreviate by $D1_{|D| \leq \alpha_p}$, and let $\tilde{\Gamma}' = \Gamma' 1_{|S| \leq \beta_p}$ (this is understood using the singular value decomposition of $\Gamma'$, where we keep the singular values that are less than $\sqrt{\beta_p}$ and replace the others by 0).

We see that $\operatorname{rank}(D - \tilde{D}) = \sum_{i=1}^n 1_{\lambda_i^2 > \alpha_p}$ and, similarly, $0 \leq \operatorname{rank}(\Gamma' - \tilde{\Gamma}') \leq \sum_{i=1}^d 1_{|s_i| > \beta_p}$. Since we assumed that $G_d$ converges weakly a.s. to $G$ and $\nu_n$ converges weakly a.s. to $\nu$, we conclude that for $\alpha_p = \beta_p = \log p$, $\operatorname{rank}(\Gamma' - \tilde{\Gamma}')/p \to 0$ a.s. and $\operatorname{rank}(D - \tilde{D})/p \to 0$ a.s. Here, it is important that $d/p$ and $p/n$ have finite nonzero limits.

So, to prove the theorem, it is sufficient to prove it for $D$ and $S$ bounded in operator norm by, for instance, $\log p$ since we just showed that by truncating $S$ and $D$ at these levels, we will not change the limiting spectral distribution of the matrices of interest, provided it exists.

### 3.3.2. *Proof of Theorem 2.*
As explained in Section 3.3.1, we can, and do, assume that all of the eigenvalues of $S = \Gamma\Gamma'$ are less than $\log p$ and, similarly, we assume that $|\lambda_i| < \sqrt{\log p}$ since, as we have explained, these assumptions do not affect the limiting spectral distribution of $B_n$. We also recall that we assume that $\|\|\Sigma\|\|_2 \leq \log(p)$. We call the spectral measures obtained after truncation $\tilde{G}_d$, $\tilde{H}_d$ and $\tilde{\nu}_n$, to keep track of the modifications we have induced by truncation. However, to avoid cumbersome notation, we use $S$, $T$ and $\Gamma$ to refer to the matrices we deal with. ($\tilde{S}$ and $\tilde{T}$ might have been more appropriate, but the notation would be too heavy.) The approach we use follows the "rank one perturbation" approach developed in [45] and [44].

We remind the reader that, as clearly explained in, for instance, [22], one can show vague convergence of distributions by showing pointwise convergence of Stieltjes transforms. This is the approach we take and we will therefore show convergence at fixed $z$ of the Stieltjes transforms of interest. Finally, we will need a tightness result to go from vague to weak convergence. We now turn to the actual proof.

Recall that $u_k = \sqrt{\theta_n/n}\lambda_k \Gamma r_k$, $\theta_n = d/p$ and $B_n = \sum_{i=1}^n u_i u_i'$. We define

- $B_{(k)} = B_n - u_k u_k' = \sum_{i \neq k} u_i u_i'$,

- $M_k = (B_n - u_k u_k' - z \operatorname{Id}_d)^{-1} = (B_{(k)} - z \operatorname{Id}_d)^{-1}$,
- $\mathcal{M}_n = (B_n - z \operatorname{Id}_d)^{-1}$

and

$$\beta(z) = \frac{1}{n} \sum_{k=1}^{n} \frac{\theta_n \lambda_k^2}{1 + u_k' M_k u_k}.$$

We note that $B_n$ is $d \times d$, as are all of the other matrices involved here. Using the first resolvent identity $A^{-1} - B^{-1} = A^{-1}(B - A)B^{-1}$ and the fact that (see [44])

$$(2) \qquad B_n(B_n - z\operatorname{Id}_d)^{-1} = \operatorname{Id}_d + z(B_n - z\operatorname{Id}_d)^{-1} = \sum_{k=1}^{n} \frac{u_k u_k' M_k}{1 + u_k' M_k u_k},$$

we have

$$(\beta(z)T - z\operatorname{Id}_d)^{-1} - (B_n - z\operatorname{Id}_d)^{-1}$$
$$= (\beta(z)T - z\operatorname{Id}_d)^{-1} \left[ \sum_{k=1}^{n} \frac{u_k u_k' M_k}{1 + u_k' M_k u_k} - \beta(z)T(B_n - z\operatorname{Id}_d)^{-1} \right]$$

and, hence,

$$(\beta(z)T - z\operatorname{Id}_d)^{-1} - (B_n - z\operatorname{Id}_d)^{-1}$$
$$= \sum_{k=1}^{n} \frac{1}{1 + u_k' M_k u_k} \left[ (\beta(z)T - z\operatorname{Id}_d)^{-1} u_k u_k' M_k \right.$$
$$\left. - \frac{\theta_n}{n} \lambda_k^2 (\beta(z)T - z\operatorname{Id}_d)^{-1} T(B_n - z\operatorname{Id}_d)^{-1} \right].$$

Taking traces and dividing by $d$, we get

$$\int \frac{d\widetilde{H}_d(\tau)}{\beta(z)\tau - z} - m_n(z)$$
$$(3) \qquad = \frac{1}{d} \sum_{k=1}^{n} \frac{1}{1 + u_k' M_k u_k} \left[ u_k' M_k (\beta(z)T - z\operatorname{Id}_d)^{-1} u_k \right.$$
$$\left. - \frac{\theta_n}{n} \lambda_k^2 \operatorname{trace}((\beta(z)T - z\operatorname{Id}_d)^{-1} T\mathcal{M}_n) \right].$$

Now, using, for instance, equation (2.3) in [44], we easily obtain

$$\left| \frac{1}{1 + u_k' M_k u_k} \right| \leq \frac{|z|}{v}.$$

On the other hand, it is clear that $\operatorname{Im}[\beta(z)] \leq 0$. As a matter of fact, the eigenvalues of $M_k$ all have positive imaginary part [if $z = u + iv$, they are $1/(\lambda_j(B_{(k)}) - u - iv)]$. Note, also, that $\|\|M_k\|\|_2 \le 1/v$. According to our first remark, the imaginary part of $1 + u_k' M_k u_k$ is positive and the imaginary part of $1/(1 + u_k' M_k u_k)$ is negative. Hence, the imaginary part of the eigenvalues of $\beta(z)T - z\operatorname{Id}_d$ is smaller than $-v$ ($T$ is positive semidefinite) and their modulus is greater than $v$. Therefore,

$$\|\operatorname{Re}[(\beta(z)T - z\operatorname{Id}_d)^{-1}]\|\|_2 \le \frac{1}{v} \quad \text{and} \quad \|\|\operatorname{Im}[(\beta(z)T - z\operatorname{Id}_d)^{-1}]\|\|_2 \le \frac{1}{v}.$$

Now, $\beta(z)$ depends on all the $u_k$'s in a nontrivial way, so we cannot apply our concentration results directly. Also, recall that $T$ is positive semidefinite, so we can write $T = \sum_{i=1}^d \tau_i e_i e_i'$, with $\tau_i \ge 0$. So, if $b(z)$ is another complex number, we have

$$(\beta(z)T - z\operatorname{Id}_d)^{-1} - (b(z)T - z\operatorname{Id}_d)^{-1} = \sum_{i=1}^d \frac{\tau_i(b(z) - \beta(z))}{(\tau_i b(z) - z)(\tau_i b(z) - z)} e_i e_i'$$

and

$$T^m[(\beta(z)T - z\operatorname{Id}_d)^{-1} - (b(z)T - z\operatorname{Id}_d)^{-1}]T^l$$
$$= \sum_{i=1}^d \frac{\tau_i^{l+m+1}(b(z) - \beta(z))}{(\tau_i b(z) - z)(\tau_i b(z) - z)} e_i e_i'.$$

Therefore, if $b(z)$ is such that $|\beta(z) - b(z)| \le \varepsilon$ and $\operatorname{Im}[b(z)] \le 0$, we have,

$$(4) \qquad \|\|(\beta(z)T - z\operatorname{Id}_d)^{-1} - (b(z)T - z\operatorname{Id}_d)^{-1}\|\|_2 \le \frac{\varepsilon\|\|T\|\|_2}{v^2},$$

$$(5) \quad |u_k' M_k(\beta(z)T - z\operatorname{Id}_d)^{-1} u_k - u_k' M_k(b(z)T - z\operatorname{Id}_d)^{-1} u_k| \le \frac{1}{v^3}\varepsilon\|\|T\|\|_2\|u_k\|_2^2$$

and

$$(6) \qquad \left|\frac{1}{d}\operatorname{trace}(T^l M_k[(\beta(z)T - z\operatorname{Id}_d)^{-1} - (b(z)T - z\operatorname{Id}_d)^{-1}])\right| \le \frac{4\|\|T\|\|_2^{l+1}\varepsilon}{v^3},$$

by decomposing the matrices appearing in the trace into real and imaginary parts, which are both symmetric in this instance, and using a well-known result (see, e.g. [2], Theorem A.4.7) on bounds of the trace of a product of symmetric matrices.

Consider

$$b_n(z) = \frac{\theta_n}{n}\sum_{k=1}^n \frac{\lambda_k^2}{1 + \xi_n\lambda_k^2 \mathbf{E}(\Omega_1(z))} \qquad \text{with } \Omega_1(z) = \frac{1}{d}\operatorname{trace}(T(B_n - z\operatorname{Id}_p)^{-1}).$$

Since $T$ is positive semidefinite, it is clear that $\operatorname{Im}[b_n(z)] \le 0$. Our aim in the next few lines is to show that $|b_n(z) - \beta(z)|$ is small. Recall that, according to Lemma 2.6 in [45], we have, for any Hermitian matrix $A$,

$$|\operatorname{trace}((M_k - \mathcal{M}_n)A)| \le \frac{\|\|A\|\|_2}{v}.$$

Applying Corollary 4, page 27, to $r_i$ and the random matrices $\Gamma' M_i \Gamma$, whose norms are bounded by $K = \log(p)/v$, we see that, for any fixed $\varepsilon > 0$ and $\delta > 0$,

$$\max_i \left| \frac{1}{d} r_i' \Gamma' M_i \Gamma r_i - \mathbf{E}(\Omega_1(z)) \right| \leq \varepsilon \frac{(\log p)^{(3/2 + 1/b + \delta)}}{v \sqrt{pc(p)} \theta_n} \triangleq \varepsilon \gamma_p \qquad \text{a.s.}$$

When this happens, we have, if we let $\alpha_k = r_k' \Gamma' M_k \Gamma r_k / d = u_k' M_k u_k / (\xi_n \lambda_k^2)$ and $\alpha = \mathbf{E}(\Omega_1(z))$,

$$|\beta(z) - b_n(z)| = \left| \frac{\theta_n}{n} \sum_{k=1}^n \left( \frac{\lambda_k^2}{1 + \xi_n \lambda_k^2 \alpha_k} - \frac{\lambda_k^2}{1 + \xi_n \lambda_k^2 \alpha} \right) \right|$$

$$\leq \frac{\xi_n \theta_n}{n} \sum_{k=1}^n \frac{\lambda_k^4 \varepsilon \gamma_p}{|(1 + \xi_n \lambda_k^2 \alpha_k)(1 + \xi_n \lambda_k^2 \alpha)|} \leq \frac{\xi_n \theta_n \varepsilon |z|^2 \gamma_p}{n v^2} \sum_{k=1}^n \lambda_k^4.$$

So, finally, since $|\lambda_k| \leq \sqrt{\log(p)}$, we have, for $C(z)$ independent of $p$,

$$|\beta(z) - b_n(z)| \leq C(z) \varepsilon \frac{(\log p)^{4 + 1/b}}{\sqrt{pc(p)}} \qquad \text{a.s.}$$

Therefore, since $\|\!\|T\|\!\|_2 \leq (\log p)^2$, using equation (4), we have

$$\left| \int \frac{d\widetilde{H}_d(\tau)}{\beta(z)\tau - z} - \frac{d\widetilde{H}_d(\tau)}{b_n(z)\tau - z} \right| \leq \varepsilon \frac{\xi_n \theta_n |z|^2 \gamma_p (\log p)^2}{n v^2} \sum_{k=1}^n \lambda_k^4$$

$$\leq \varepsilon \frac{\xi_n \theta_n |z|^2 \gamma_p}{v^2} (\log p)^4 \to 0 \qquad \text{a.s.}$$

Similarly, using our concentration bounds from Corollary 4 applied to

$$\left[ u_k' M_k (b_n(z) T - z \operatorname{Id}_d)^{-1} u_k / \lambda_k^2 - \frac{\theta_n}{n} \operatorname{trace}((b_n(z) T - z \operatorname{Id}_d)^{-1} T \mathcal{M}_n) \right],$$

we see that a.s.

$$\max_{1 \leq k \leq p} \left| \left[ u_k' M_k (b_n(z) T - z \operatorname{Id}_d)^{-1} u_k / \lambda_k^2 - \frac{\theta_n}{n} \operatorname{trace}((b_n(z) T - z \operatorname{Id}_d)^{-1} T \mathcal{M}_n) \right] \right|$$

$$\leq \frac{\varepsilon \xi_n \gamma_p}{v}$$

and therefore

$$|\Delta(b_n(z))| \triangleq \left| \frac{1}{d} \sum_{k=1}^n \frac{1}{1 + u_k' M_k u_k} \left[ u_k' M_k (b_n(z) T - z \operatorname{Id}_d)^{-1} u_k \right. \right.$$

(7)
$$\left. \left. - \frac{\theta_n}{n} \lambda_k^2 \operatorname{trace}((b_n(z) T - z \operatorname{Id}_d)^{-1} T \mathcal{M}_n) \right] \right|$$

$$\leq C \varepsilon \gamma_p \frac{|z|}{v^2} \frac{\xi_n}{d} \sum_{k=1}^n \lambda_k^2 \to 0 \qquad \text{a.s.}$$

We now need to show that $|\Delta(\beta(z))|$ tends to 0 almost surely. To do so, we study $|\Delta(\beta(z)) - \Delta(b_n(z))|$. Using equation (6), we have

$$\left| \frac{1}{d} \sum_{k=1}^{n} \frac{\lambda_k^2}{1 + u_k' M_k u_k} \left[ \frac{\theta_n}{n} \operatorname{trace}((\beta(z)T - z \operatorname{Id}_d)^{-1} T \mathcal{M}_n) \right. \right.$$
$$\left. \left. - \frac{\theta_n}{n} \operatorname{trace}((b_n(z)T - z \operatorname{Id}_d)^{-1} T \mathcal{M}_n) \right] \right|$$
$$\leq \theta_n \frac{4 \|\|T\|\|_2^2 |b_n(z) - \beta(z)| |z|}{v^4} \frac{1}{n} \sum_{k=1}^{n} \lambda_k^2 \to 0 \qquad \text{a.s.}$$

On the other hand, using equation (5), we find that

$$\left| \frac{1}{d} \sum_{k=1}^{n} \frac{1}{1 + u_k' M_k u_k} [u_k' M_k (b_n(z)T - z \operatorname{Id}_d)^{-1} u_k - u_k' M_k (\beta(z)T - z \operatorname{Id}_d)^{-1} u_k] \right|$$
$$\leq \frac{|b_n(z) - \beta(z)| |z|}{v^4} \|\|T\|\|_2 \frac{1}{n} \sum_{k=1}^{n} \lambda_k^2 \frac{1}{p} r_k' \Gamma' \Gamma r_k.$$

Applying Lemma 8, with the matrix $M = \Gamma' \Gamma$, whose operator norm is bounded by $\log(p)$, we get, as above, that

$$\max_k \left| \frac{1}{p} r_k' \Gamma' \Gamma r_k - \frac{1}{p} \operatorname{trace}(T) \right| \leq \varepsilon \frac{(\log p)^{3/2 + 1/b + \delta}}{\sqrt{pc(p)}} \qquad \text{a.s.}$$

Because $\|\|T\|\|_2 \leq (\log p)^2$, we conclude that a.s.

$$\max_k \left| \frac{1}{p} r_k' \Gamma \Gamma' r_k \right| \leq 3 (\log p)^2.$$

Therefore,

$$\frac{|b_n(z) - \beta(z)| |z|}{v^4} \|\|T\|\|_2 \frac{1}{n} \sum_{k=1}^{n} \lambda_k^2 \frac{1}{p} r_k' \Gamma' \Gamma r_k \to 0 \qquad \text{a.s.}$$

Hence,

$$|\Delta(\beta(z))| \leq |\Delta(\beta(z)) - \Delta(b_n(z))| + |\Delta(b_n(z))| \to 0 \qquad \text{a.s.}$$

Since, by equation (3),

$$\Delta(\beta(z)) = \int \frac{d\widetilde{H}_d(\tau)}{\beta(z)\tau - z} - m_n(z),$$

we can finally conclude that

$$\int \frac{d\widetilde{H}_d(\tau)}{b_n(z)\tau - z} - m_n(z) \to 0 \qquad \text{a.s.}$$

This corresponds to the first part of the theorem. Now, note that $\mathrm{Im}[b_n(z)] \leq 0$ and therefore $|1/(b_n(z)\tau - z)| \leq 1/v$. Because $\int |d\widetilde{H}_d(\tau) - dH_d(\tau)| \to 0$, we conclude that

$$\int \frac{dH_d(\tau)}{b_n(z)\tau - z} - m_n(z) \to 0 \qquad \text{a.s.}$$

To get to the second part of the theorem, we instead consider

$$T(\beta(z)T - z\,\mathrm{Id}_d)^{-1} - T(B_n - z\,\mathrm{Id}_d)^{-1}.$$

Taking traces and dividing by $d$, we get

$$\int \frac{\tau\, d\widetilde{H}_d(\tau)}{\tau\beta(z) - z} - \frac{1}{d}\,\mathrm{trace}(T(B_n - z\,\mathrm{Id}_d)^{-1}).$$

To control this quantity, we can use the same expansions we used before, everywhere replacing $(\beta(z)T - z\,\mathrm{Id}_d)^{-1}$ by $T(\beta(z)T - z\,\mathrm{Id}_d)^{-1}$. This has the effect of multiplying the upper bounds by $\|\|T\|\|_2$, which, under our assumptions, is bounded by $(\log p)^2$. So, we conclude that

$$\int \frac{\tau\, d\widetilde{H}_d(\tau)}{\tau b_n(z) - z} - \Omega_1(z) \to 0 \qquad \text{a.s.}$$

Now, the result we obtained using Azuma's inequality shows clearly (see Remark 1) that

$$\Omega_1(z) - \mathbf{E}(\Omega_1(z)) \to 0 \qquad \text{a.s.}$$

Letting $w_n(z) = \mathbf{E}(\Omega_1(z))$, we have shown that

$$\begin{cases} \displaystyle\int \frac{\tau\, d\widetilde{H}_d(\tau)}{\tau \int \theta_n \lambda^2\, d\bar{\nu}_n(\lambda)/(1 + \xi_n\lambda^2 w_n(z)) - z} - w_n(z) \to 0 \qquad \text{a.s.} \quad \text{and} \\[4mm] \displaystyle\int \frac{d\widetilde{H}_d(\tau)}{\tau \int \theta_n \lambda^2\, d\bar{\nu}_n(\lambda)/(1 + \xi_n\lambda^2 w_n(z)) - z} - m_n(z) \to 0 \qquad \text{a.s.} \end{cases}$$

(8)

• *Subsequence argument to reach the conclusion of Theorem 2*

We now need to turn to technical arguments to get from the statement of equation 8 to that of Theorem 2. Because of our assumption that $\int \tau\, dH_d(\tau) < K$ for all $d$ (or $p$, which is equivalent), with $K$ fixed and independent of $d$, we see that $|w_n(z)| \leq \mathrm{trace}(T)/(dv) < K/v$. So, at fixed $z$, $w_n(z)$ is bounded. From this sequence, let us extract a convergent subsequence $w_{m(n)}(z)$, or $w_m$ for short, that converges to $w$. Through tightness arguments (see below), we see that $w \in \mathbb{C}^+$. We will now show that $w(z)$ satisfies

$$\int \frac{\tau\, dH(\tau)}{\tau \int \theta\lambda^2\, d\nu(\lambda)/(1 + \xi\lambda^2 w(z)) - z} - w(z) = 0$$

and that there is a unique solution to this equation in $\mathbb{C}^+$. Let $b_m(z) = \int \theta_m \times \lambda^2 \, d\tilde{\nu}_m(\lambda)/(1 + \xi_m \lambda^2 w_m(z))$. We first show that $b_m \to b = \int \theta \lambda^2 \, d\nu(\lambda)/(1 + \xi \lambda^2 w(z))$. To do so, note that $\lambda^2/(1 + w_m \lambda^2) - \lambda^2/(1 + w\lambda^2) = (w - w_m)\lambda^4/[(1 + w\lambda^2)(1 + w_m\lambda^2)]$. Now, because $w_m \to w \in \mathbb{C}^+$, their imaginary parts are uniformly bounded below by $\delta$, from which we conclude that, if $w_m \to w \in \mathbb{C}^+$,

$$\int \frac{\lambda^2 \, d\tilde{\nu}_m}{1 + w_m \lambda^2} - \int \frac{\lambda^2 \, d\tilde{\nu}_m}{1 + w\lambda^2} \to 0.$$

On the other hand, for $w \in \mathbb{C}^+$, $\lambda^2/(1 + w\lambda^2)$ is a bounded continuous function of $\lambda$. Since $\nu_m \Rightarrow \nu$ and, therefore, $\tilde{\nu}_m \Rightarrow \nu$, we conclude that

$$\int \frac{\lambda^2 \, d\tilde{\nu}_m}{1 + a\lambda^2} \to \int \frac{\lambda^2 \, d\nu}{1 + a\lambda^2}.$$

Therefore, since $\theta_m \to \theta$, $b_m(z) \to b(z)$. Because we have assumed that $\nu \neq 0$, we have $b(z) \in \mathbb{C}^-$. By essentially the same arguments, using the fact that $|\operatorname{Im}[b_m(z)]|$ is bounded below by $\delta$ and $b(z) \in \mathbb{C}^-$, we conclude that

$$\int \frac{\tau \, d\tilde{H}_{d(m(n))}(\tau)}{\tau b_{m(n)}(z) - z} - \int \frac{\tau \, dH(\tau)}{\tau b(z) - z} \to 0.$$

In other words,

$$\int \frac{\tau \, dH(\tau)}{\tau b(z) - z} - w(z) = 0,$$

where

$$b(z) = \int \frac{\theta \lambda^2 \, d\nu(\lambda)}{1 + \xi \lambda^2 w(z)}.$$

Similarly, we can show that along this subsequence,

$$\int \frac{dH_d(\tau)}{\tau b_m(z) - z} \to \int \frac{dH(\tau)}{\tau b(z) - z}$$

and so we also get the first equation in Theorem 2.

• *Uniqueness of possible limit*

We now prove that there is a unique solution in $\mathbb{C}^+$ to the equation characterizing $w$, the only question remaining to tackle being uniqueness. To do so, we employ an argument similar to that given in [45], although the details are slightly different.

Suppose we have two solutions in $\mathbb{C}^+$ to the equation characterizing $w(z)$. Let us call them $w_1$ and $w_2$, where $b_1$ and $b_2$ are the corresponding $b$'s. We have

$$w_1 - w_2 = \int \left( \frac{\tau}{\tau b_1 - z} - \frac{\tau}{\tau b_2 - z} \right) dH(\tau)$$

$$= (b_2 - b_1) \int \frac{\tau^2}{(\tau b_1 - z)(\tau b_2 - z)} \, dH(\tau)$$

$$= \theta(w_1 - w_2) \int \frac{\lambda^4 \xi \, d\nu(\lambda)}{(1 + \xi \lambda^2 w_1(z))(1 + \xi \lambda^2 w_2(z))}$$

$$\times \int \frac{\tau^2}{(\tau b_1 - z)(\tau b_2 - z)} \, dH(\tau).$$

Let us call $f$ the quantity multiplying $w_1 - w_2$ in the previous equation. We want to show that $|f| < 1$. As in [45], using Hölder's inequality, we have, given that $\theta > 0$,

$$|f| \leq \left( \theta \int \frac{\lambda^4 \xi \, d\nu(\lambda)}{|1 + \xi \lambda^2 w_1(z)|^2} \int \frac{\tau^2}{|\tau b_1 - z|^2} \, dH(\tau) \right)^{1/2}$$

$$\times \left( \theta \int \frac{\lambda^4 \xi \, d\nu(\lambda)}{|1 + \xi \lambda^2 w_2(z)|^2} \int \frac{\tau^2}{|\tau b_2 - z|^2} \, dH(\tau) \right)^{1/2}.$$

Let us write $w_1 = a + ic$, $z = u + iv$ and $b_1 = \alpha - i\gamma$. By writing the definition of $b_1$ in terms of $w_1$, we see immediately that

$$\gamma = c \int \frac{\theta \xi \lambda^4}{|1 + \xi \lambda^2 w_1|^2} \, d\nu(\lambda),$$

so $\int \frac{\theta \xi \lambda^4}{|1 + \xi \lambda^2 w_1|^2} \, d\nu(\lambda) = -\operatorname{Im}[b_1] / \operatorname{Im}[w_1]$. Since $\nu \neq 0$ by our assumptions, we see that $\gamma > 0$. On the other hand, using the definition of $w_1$ in terms of $b_1$, we see that

$$\operatorname{Im}[w_1] = \int -\operatorname{Im}[b_1] \frac{\tau^2}{|\tau b_1 - z|^2} \, dH(\tau) + \operatorname{Im}[z] \int \frac{\tau}{|\tau b_1 - z|^2} \, dH(\tau)$$

and, therefore, $\operatorname{Im}[w_1] > -\operatorname{Im}[b_1] \int \frac{\tau^2}{|\tau b_1 - z|^2} \, dH(\tau)$ since $H \neq 0$.

Hence,

$$\left( \int \frac{\theta \lambda^4 \xi \, d\nu(\lambda)}{|1 + \xi \lambda^2 w_1(z)|^2} \int \frac{\tau^2}{|\tau b_1 - z|^2} \, dH(\tau) \right)^{1/2} < 1$$

and $|f| < 1$. We conclude that $w_2 = w_1$, so there is at most one solution to the equation characterizing $w$.

• *Tightness of $B_n$ and consequences for $w$*

Finally, we need to show that the spectral distribution $F^{B_n}$ is tight a.s. and deduce consequences for $w$. It is shown (via Lemma 2.3) in [45] that if $B_n = T_n^{1/2} Y_n^* Y_n T_n^{1/2}$, if the spectral distributions of the $T_n$'s form a tight sequence and so do the spectral distributions of the $Y_n^* Y_n$'s, then $F^{B_n}$ form a tight sequence. We note that in our case $B_n = \Gamma R_n^* D_n^2 R_n \Gamma' / n$, which, up to a number of zeros, has the same eigenvalues as $S^{1/2} R_n^* D_n^2 R_n S^{1/2} / n$;

we temporarily denote by $R_n$ the matrix containing our vectors $r_i$. So, all we have to show is that $F^{R_n^* D_n^2 R_n / n}$ forms a tight sequence. Note that our assumption on the convergence of the spectral distribution of the $\lambda$'s implies that the spectral distributions of the $D_n^2$'s form a tight sequence. So, all we need do in order to conclude is to show that $F^{R_n^* R_n / n}$ also forms a tight sequence. But we showed this in Corollary 5. So, $F^{B_n}$ forms a tight sequence, a.s. Recall that when $\mathrm{trace}(\Sigma)/p$ is uniformly bounded by $K$, we showed in Corollary 5 that a.s. $F^{R_n^* R_n / n}([M, \infty)) \leq (K+1)/M$. So, for any $\varepsilon$, we can find $M_\varepsilon$ such that $F^{B_n}[M_\varepsilon, \infty) < \varepsilon$, a.s. Using the second inequality in Lemma 2.3 in [45] and the fact that $H$ and $\nu$ are deterministic, as well as the fact that if $X_n \Rightarrow X$ and $C$ is a closed set, $\limsup P(X_n \in C) \leq P(X \in C)$, we see that $M_\varepsilon$ can be chosen uniformly in $\omega$.

We now want to show that $w \in \mathbb{C}^+$; to do so, we will show that a.s., $\mathrm{Im}[w_n]$ is bounded away from zero. Note that $\mathrm{Im}[(B_n - z \,\mathrm{Id})^{-1}]$ is a symmetric matrix. Its eigenvalues, which we denote by $a_k$, are, if $l_k$ denote the eigenvalues of $B_n$, $v/((l_k - u)^2 + v^2) \geq v/(2(l_k^2 + u^2) + v^2)$. Assume that $a_1 \geq a_2 \geq \cdots \geq a_d$. Using Theorem A.4.7 in [2], we see that, if we denote by $\tau_i$ the decreasingly ordered eigenvalues of $T$,

$$\mathrm{Im}[\Omega_1(z)] = \mathrm{Im}\left[\frac{1}{d}\,\mathrm{trace}(T(B_n - z \,\mathrm{Id}_p)^{-1})\right] \geq \frac{1}{d}\sum_{i=1}^{d}\tau_i a_{d-i+1}.$$

Now, all we need to show is that a.s. a fixed nonzero proportion of $\tau_i a_{d-i}$ stay bounded away from 0. Because $H \neq 0$, we can find $\eta$ such that $H(\eta, \infty) > \varepsilon$ for some $\varepsilon > 0$. Let us choose such an $\varepsilon \neq 0$. In particular, the proportion of indices for which $\tau_i > \eta$ is a.s. greater than $\varepsilon$ because $\liminf H_d(\eta, \infty) \geq H(\eta, \infty)$, a.s. For this $\varepsilon$, we can find $m_\varepsilon < \infty$ such that $F^{B_n}[0, m_\varepsilon] \geq 1 - \varepsilon/2$ a.s., from our arguments above. So, the proportion of $i$'s such that $a_{d-i+1} \geq v/(2(m_\varepsilon^2 + u^2) + v^2)$ is greater than $1 - \varepsilon/2$. So, the proportion of $i$'s for which both $\tau_i > \eta$ and $a_{d-i+1} \geq v/(2(m_\varepsilon^2 + u^2) + v^2)$ must be greater than $\varepsilon/2$, a.s. Hence, $\mathrm{Im}[\Omega_1(z)] \geq \delta > 0$, a.s. Now, we saw that $w_n(z) = \mathbf{E}(\Omega_1(z))$ is such that $w_n(z) - \Omega_1(z) \to 0$ a.s. Hence, $\mathrm{Im}[w_n(z)] \geq \delta > 0$ a.s. and we can conclude that $w$ is in $\mathbb{C}^+$. So, the subsequence argument given above is valid and we can go from the result of equation (8) to the main result of Theorem 2.

So, using the connection between pointwise convergence of Stieltjes transforms (see [22]) and vague convergence, we have shown that the spectral distribution of $B_n$ converges vaguely a.s. to a nonrandom distribution, which is uniquely characterized by the system of equations described in Theorem 2. The a.s. tightness result we obtained for $F^{B_n}$ ensures that the limiting spectral distribution of $B_n$ is a probability measure and, hence, we have a.s. weak convergence, as announced in Theorem 2. This completes the proof.

**4. Conclusion.** We have shown that the concentration of measure phenomenon can be seen as an essential tool in the understanding of the behavior of the limiting spectral distributions of a number of random matrix models.

Motivated by applications, we have used one aspect of this phenomenon to deduce spectral properties of sample correlation matrices from the corresponding properties for sample covariance matrices. On the other hand, for more complicated models, we have generalized known results about random covariance-type matrices to sample covariance matrices computed from elliptically distributed data, a type of assumption that is popular in financial modeling and, further, to generalized elliptically distributed data. We have done this almost entirely from concentration properties of certain quadratic forms. An interesting aspect of the proof is that it leads to new results for data coming from distributions for which the dependence between entries of the data vector cannot be broken up in a linear fashion. The concentration approach also highlights the fact that data vectors coming from a distribution having the dimension-free concentration property we used repeatedly have, after proper normalization, almost the same norm and are almost orthogonal to one another (in concrete terms, this remark applies to models considered in Theorems 1 and 2 when $\lambda_i = 1$). Since this peculiar geometric feature may not be present in data sets to be analyzed, practitioners should probably perform corresponding diagnostic checks before relying on random matrix results of the type discussed in this and other papers.

Interestingly, in all of the models considered, the results tell us that only the covariance or the correlation between the entries of the data vector matters and the more complicated dependence structure is irrelevant as far as limiting distributions of eigenvalues are concerned.

## APPENDIX

**On the covariance matrix of data distributed according to a Gaussian copula.** The concentration results we developed in Section 3.2 require the covariance matrices of the data at stake to be bounded in operator norm by $\log(p)$. To be able to apply the results to data distributed according to a Gaussian copula, we therefore need to show that this is the case. We have the following fact.

FACT 2. *Suppose $r$ is distributed according to a Gaussian copula with corresponding correlation matrix $R$. Then, if $\Sigma$ is the covariance matrix of $r$, we have*

$$|||\Sigma|||_2 \leq \frac{1}{2\pi}(|||R|||_2/2 + 4|||R|||_2^2(\pi/6 - 1/2)).$$

Proof. Recall that if $r$ is distributed according to a Gaussian copula with corresponding correlation matrix $R$, $r$ can be generated in the following way: draw $v$ according to a multivariate normal $\mathcal{N}(0, R)$. Because $R$ is a correlation matrix, $v_i$, the $i$th entry of $v$ is $\mathcal{N}(0, 1)$. Now, calling $\Phi$ the c.d.f. of the standard normal distribution, $r_i = \Phi(v_i)$.

We also recall the standard fact (see [37], Definition 5.28, Proposition 5.29 and Theorem 5.36) that

$$\Sigma(i, j) = \mathrm{cov}(r_i, r_j) = \frac{1}{2\pi} \arcsin(R_{i,j}/2).$$

Note that $|R_{i,j}| \leq 1$. Recall the series expansion of $\arcsin(x)$, valid for $x < 1$:

$$\arcsin(x) = \sum_{n=0}^{\infty} u_n x^{2n+1}, \qquad \text{with } u_n = \frac{(2n)!}{4^n (n!)^2 (2n+1)}.$$

Denote by $E \circ M$ the Hadamard (i.e. entrywise) product of matrices $E$ and $M$. In [15], it is shown that if $M$ has nonnegative entries, and $M$ and $E$ are symmetric, then

$$\|\!|E \circ M|\!\|_2 \leq \max_{i,j}(|E_{i,j}|)\|\!|M|\!\|_2.$$

Call $g(R/2)$ the matrix with entries $\arcsin(R_{i,j}/2)$. Then

$$g(R/2) = \frac{R}{2} + \sum_{n=1}^{\infty} u_n \left[\frac{R}{2}\right]^{\circ 2n-1} \circ \left[\frac{R}{2}\right]^{\circ 2},$$

where $A^{\circ n}$ is the $n$th Hadamard power (i.e. entrywise) of matrix $A$. Now $\max_{i,j} |R_{i,j}/2| \leq 1/2$, so

$$\max_{i,j} |(R_{i,j}/2)^{2n-1}| \leq \frac{1}{2^{2n-1}}.$$

Now, using e.g. [7], Problem I.6.13, we have

$$\|\!|R \circ R|\!\|_2 \leq \|\!|R|\!\|_2^2.$$

We therefore have

$$\left\|\!\left|g\left(\frac{R}{2}\right)\right|\!\right\|_2 \leq \frac{\|\!|R|\!\|_2}{2} + \sum_{n=1}^{\infty} u_n \frac{1}{2^{2n-1}} \|\!|R|\!\|_2^2$$

$$\leq \frac{\|\!|R|\!\|_2}{2} + 4\|\!|R|\!\|_2^2 \sum_{n=1}^{\infty} u_n \frac{1}{2^{2n+1}}$$

$$\leq \frac{\|\!|R|\!\|_2}{2} + 4\|\!|R|\!\|_2^2 (\arcsin(1/2) - 1/2).$$

The result follows since $\arcsin(1/2) = \pi/6$. $\quad\square$

**Acknowledgments.** Hospitality and support from SAMSI in the fall of 2006 are gratefully acknowledged. I would like to thank James Mingo and Jack Silverstein for an interesting conversation while there. I would also like to thank an anonymous referee and the Editor for constructive comments that led to a substantial improvement of the paper.

## REFERENCES

[1] Anderson, G. W. and Zeitouni, O. (2006). A CLT for a band matrix model. *Probab. Theory Related Fields* **134** 283–338. MR2222385

[2] Anderson, T. W. (2003). *An Introduction to Multivariate Statistical Analysis*, 3rd ed. *Wiley Series in Probability and Statistics.* Wiley, Hoboken, NJ. MR1990662

[3] Bai, Z. D. (1999). Methodologies in spectral analysis of large-dimensional random matrices, a review. *Statist. Sinica* **9** 611–677. MR1711663

[4] Bai, Z. D. and Silverstein, J. W. (1998). No eigenvalues outside the support of the limiting spectral distribution of large-dimensional sample covariance matrices. *Ann. Probab.* **26** 316–345. MR1617051

[5] Bai, Z. D. and Silverstein, J. W. (2004). CLT for linear spectral statistics of large-dimensional sample covariance matrices. *Ann. Probab.* **32** 553–605. MR2040792

[6] Baik, J., Ben Arous, G. and Péché, S. (2005). Phase transition of the largest eigenvalue for nonnull complex sample covariance matrices. *Ann. Probab.* **33** 1643–1697. MR2165575

[7] Bhatia, R. (1997). *Matrix Analysis. Graduate Texts in Mathematics* **169**. Springer, New York. MR1477662

[8] Bickel, P. J. and Levina, E. (2008). Regularized estimation of large covariance matrices. *Ann. Statist.* **36** 199–227. MR2387969

[9] Boutet de Monvel, A., Khorunzhy, A. and Vasilchuk, V. (1996). Limiting eigenvalue distribution of random matrices with correlated entries. *Markov Process. Related Fields* **2** 607–636. MR1431189

[10] Boyd, S. and Vandenberghe, L. (2004). *Convex Optimization.* Cambridge Univ. Press, Cambridge. MR2061575

[11] Burda, Z., Görlich, A., Jarosz, A. and Jurkiewicz, J. (2004). Signal and noise in correlation matrix. *Phys. A* **343** 295–310. MR2094415

[12] Burda, Z., Jurkiewicz, J. and Wacław, B. (2005). Spectral moments of correlated Wishart matrices. *Phys. Rev. E* **71**.

[13] Campbell, J., Lo, A. and MacKinlay, C. (1996). *The Econometrics of Financial Markets.* Princeton Univ. Press, Princeton, NJ.

[14] El Karoui, N. (2003). On the largest eigenvalue of Wishart matrices with identity covariance when $n, p$ and $p/n \to \infty$. Available at arXiv:math.ST/0309355.

[15] El Karoui, N. (2007a). The spectrum of kernel random matrices. Technical Report 748, Dept. Statistics, UC Berkeley. *Ann. Statist.* To appear.

[16] El Karoui, N. (2007). Tracy–Widom limit for the largest eigenvalue of a large class of complex sample covariance matrices. *Ann. Probab.* **35** 663–714. MR2308592

[17] El Karoui, N. (2009). Spectrum estimation for large dimensional covariance matrices using random matrix theory. *Ann. Statist.* To appear. Available at arXiv:math.ST/0609418.

[18] Fang, K. T., Kotz, S. and Ng, K. W. (1990). *Symmetric Multivariate and Related Distributions. Monographs on Statistics and Applied Probability* **36**. Chapman & Hall, London. MR1071174

[19] FORRESTER, P. J. (1993). The spectrum edge of random matrix ensembles. *Nuclear Phys. B* **402** 709–728. MR1236195

[20] FRAHM, G. and JAEKEL, U. (2005). Random matrix theory and robust covariance matrix estimation for financial data. Available at arXiv:physics/0503007.

[21] GEMAN, S. (1980). A limit theorem for the norm of random matrices. *Ann. Probab.* **8** 252–261. MR566592

[22] GERONIMO, J. S. and HILL, T. P. (2003). Necessary and sufficient condition that the limit of Stieltjes transforms is a Stieltjes transform. *J. Approx. Theory* **121** 54–60. MR1962995

[23] GIRKO, V. L. (1990). *Theory of Random Determinants. Mathematics and Its Applications (Soviet Series)* **45**. Kluwer Academic, Dordrecht. MR1080966

[24] GRAY, R. M. (2002). Toeplitz and circulant matrices: A review. Available at http://ee.stanford.edu/~gray/toeplitz.pdf.

[25] GRENANDER, U. and SZEGÖ, G. (1958). *Toeplitz Forms and Their Applications.* Univ. California Press, Berkeley. MR0094840

[26] GUIONNET, A. and ZEITOUNI, O. (2000). Concentration of the spectral measure for large matrices. *Electron. Comm. Probab.* **5** 119–136 (electronic). MR1781846

[27] JIANG, T. (2004). The limiting distributions of eigenvalues of sample correlation matrices. *Sankhyā* **66** 35–48. MR2082906

[28] JOHANSSON, K. (2000). Shape fluctuations and random matrices. *Comm. Math. Phys.* **209** 437–476. MR1737991

[29] JOHNSTONE, I. M. (2001). On the distribution of the largest eigenvalue in principal components analysis. *Ann. Statist.* **29** 295–327. MR1863961

[30] JONSSON, D. (1982). Some limit theorems for the eigenvalues of a sample covariance matrix. *J. Multivariate Anal.* **12** 1–38. MR650926

[31] LALOUX, L., CIZEAU, P., BOUCHAUD, J.-P. and POTTERS, M. (1999). Random matrix theory and financial correlations. *Int. J. Theor. Appl. Finance* **3** 391–397.

[32] LEDOIT, O. and WOLF, M. (2004). A well-conditioned estimator for large-dimensional covariance matrices. *J. Multivariate Anal.* **88** 365–411. MR2026339

[33] LEDOUX, M. (2001). *The Concentration of Measure Phenomenon. Mathematical Surveys and Monographs* **89**. Amer. Math. Soc., Providence, RI. MR1849347

[34] LI, L., TULINO, A. M. and VERDÚ, S. (2004). Design of reduced-rank MMSE multiuser detectors using random matrix methods. *IEEE Trans. Inform. Theory* **50** 986–1008. MR2094863

[35] MARČENKO, V. A. and PASTUR, L. A. (1967). Distribution of eigenvalues in certain sets of random matrices. *Mat. Sb. (N.S.)* **72** 507–536. MR0208649

[36] MARDIA, K. V., KENT, J. T. and BIBBY, J. M. (1979). *Multivariate Analysis.* Academic Press, London. MR560319

[37] MCNEIL, A. J., FREY, R. and EMBRECHTS, P. (2005). *Quantitative Risk Management: Concepts, Techniques and Tools.* Princeton Univ. Press, Princeton, NJ. MR2175089

[38] NELSEN, R. B. (2006). *An Introduction to Copulas*, 2nd ed. Springer, New York. MR2197664

[39] PAJOR, A. and PASTUR, L. (2007). On the limiting empirical measure of the sum of rank one matrices with log-concave distribution. Available at http://www.arxiv.org/abs/0710.1346.

[40] PAUL, D. (2007). Asymptotics of sample eigenstructure for a large dimensional spiked covariance model. *Statist. Sinica* **17** 1617–1642. MR2399865

[41] PAUL, D. and SILVERSTEIN, J. (2007). No eigenvalues outside the support of the limiting empirical spectral distribution of a separable covariance matrix. Available at http://www4.ncsu.edu/~jack/pub.html.

[42] RAO, N. R., MINGO, J., SPEICHER, R. and EDELMAN, A. (2007). Statistical eigeninference from large Wishart matrices. Available at arXiv:math/0701314.

[43] SCHECHTMAN, G. and ZINN, J. (2000). Concentration on the $l_p^n$ ball. In *Geometric Aspects of Functional Analysis. Lecture Notes in Math.* **1745** 245–256. Springer, Berlin. MR1796723

[44] SILVERSTEIN, J. W. (1995). Strong convergence of the empirical distribution of eigenvalues of large-dimensional random matrices. *J. Multivariate Anal.* **55** 331–339. MR1370408

[45] SILVERSTEIN, J. W. and BAI, Z. D. (1995). On the empirical distribution of eigenvalues of a class of large-dimensional random matrices. *J. Multivariate Anal.* **54** 175–192. MR1345534

[46] TALAGRAND, M. (1995). Concentration of measure and isoperimetric inequalities in product spaces. *Inst. Hautes Études Sci. Publ. Math.* **81** 73–205. MR1361756

[47] TULINO, A. and VERDÚ, S. (2004). *Random Matrix Theory and Wireless Communications. Foundations and Trends in Communications and Information Theory* **1**. Now Publishers, Boston.

[48] VAN DER VAART, A. W. (1998). *Asymptotic Statistics. Cambridge Series in Statistical and Probabilistic Mathematics* **3**. Cambridge Univ. Press, Cambridge. MR1652247

[49] VOICULESCU, D. (2000). Lectures on free probability theory. In *Lectures on Probability Theory and Statistics (Saint-Flour, 1998). Lecture Notes in Math.* **1738** 279–349. Springer, Berlin. MR1775641

[50] YIN, Y. Q., BAI, Z. D. and KRISHNAIAH, P. R. (1988). On the limit of the largest eigenvalue of the large-dimensional sample covariance matrix. *Probab. Theory Related Fields* **78** 509–521. MR950344

[51] ZHANG, L. (2006). Spectral analysis of large dimensional random matrices. Ph.D. thesis, National Univ. Singapore.

367 EVANS HALL
UC-BERKELEY
BERKELEY, CALIFORNIA 94720-3860
USA
E-MAIL: nkaroui@stat.berkeley.edu